\input amstex
\documentstyle{amsppt}
\magnification=\magstep1
\NoBlackBoxes

\pageheight{9.0truein}
\pagewidth{6.5truein}

\input xy
\xyoption{matrix}\xyoption{arrow}\xyoption{curve}
\def\arrow{\ar}

\long\def\ignore#1{#1}

\def\edge{\arrow@{-}}

\def\uloopr#1{\arrow@'{@+{[0,0]+(-4,5)} @+{[0,0]+(0,10)} @+{[0,0]+(4,5)}}
  ^{#1}}
\def\dloopr#1{\arrow@'{@+{[0,0]+(-4,-5)} @+{[0,0]+(0,-10)} @+{[0,0]+(4,-5)}}
  _{#1}}
\def\rloopd#1{\arrow@'{@+{[0,0]+(5,4)} @+{[0,0]+(10,0)} @+{[0,0]+(5,-4)}}
  ^{#1}}
\def\lloopd#1{\arrow@'{@+{[0,0]+(-5,4)} @+{[0,0]+(-10,0)} @+{[0,0]+(-5,-4)}}
  _{#1}}
\def\udotloopr#1{\arrow@{{}.>} @'{@+{[0,0]+(-4,5)} @+{[0,0]+(0,10)}
  @+{[0,0]+(4,5)}} ^{#1}}

\def\la{\Lambda}
\def\gam{\Gamma}
\def\kgam{K\Gamma}
\def\lamod{\Lambda\operatorname{-mod}}
\def\soc{\operatorname{Soc}}

\def\that{{\widehat t}}
\def\That{{\widehat T}}
\def\len{\operatorname{length}}
\def\pinf{{\Cal P}^\infty}
\def\pdim{\operatorname{p\,dim}}
\def\sup{\operatorname{sup}}
\def\S{{\Cal S}}
\def\sinf{{\Cal S}^\infty}

\def\lfd{\operatorname{l\,fin\,dim}}
\def\opio{{\Cal O}/\pi{\Cal O}}
\def\O{{\Cal O}}
\def\A{{\frak A}}
\def\B{{\frak B}}
\def\C{{\frak C}}
\def\D{{\frak D}}
\def\N{{\Bbb N}}
\def\gldim{\operatorname{gl\,dim}}
\def\Hom{\operatorname{Hom}}

\def\xbar{{\overline x}}

\topmatter

\title Analyzing the Structure of Representations via Approximations\endtitle
\rightheadtext{ANALYZING REPRESENTATIONS VIA APPROXIMATIONS}
\author Birge Zimmermann Huisgen\endauthor
\address Department of Mathematics, University of California, Santa Barbara,
CA 93106, USA\endaddress
\email birge\@math.ucsb.edu\endemail

\dedicatory Dedicated to the memory of Maurice Auslander.\enddedicatory

\thanks This work was partially supported by an N.S.F. grant.\endgraf
This paper is in final form and no version of it will be submitted
for publication elsewhere.\endthanks

\abstract Primarily this paper presents an expository report on alternatives to
the traditional methods of classifying representations of finite dimensional
algebras. Some new results illustrating such alternatives for algebras with
only finitely many isomorphism types of uniserial modules are included.
\endabstract


\endtopmatter

\document

\head 1. Introduction and notation \endhead

We discuss several avenues of approach to the structure of representations of
a finite dimensional algebra $\la$ by comparing the
target objects with the objects from a more thoroughly understood
reference class.  Our starting point is the concept of a `right
$\A$-ap\-prox\-i\-ma\-tion' of a left
$\la$-module relative to a subcategory $\A$ of $\lamod$, as introduced by
Auslander and Smal\o\ in \cite{3} (although originally not in this
terminology), as well as the subsequent work on the subject by Auslander and
Reiten
\cite{2}.

Due to \cite{2}, the following holds for any resolving
contravariantly finite subcategory $\A$ of $\lamod$: 
If there are $n$ simple left $\la$-modules, up to isomorphism, and if $A_1,
\dots, A_n$ are their minimal right $\A$-ap\-prox\-i\-ma\-tions, then a left
$\la$-module belongs to $\A$ if and only if it is isomorphic to a direct
summand of a module $M$ having a filtration $M = M_0 \supseteq M_1 \supseteq
\dots \supseteq M_l = 0$ with $M_i/M_{i+1} \in \{A_1, \dots ,A_n\}$.  In
other words, approximating the simple modules by objects in $\A$ yields
basic structural information about arbitrary objects in $\A$.  To give
a flavor of the usefulness of this type of information, e.g., towards a
homological understanding of $\A$, we point out that, in the described
situation, we have $\sup\{\pdim A \mid A \in \A\} = \sup\{\pdim A_i \mid
1\le i \le n\}$.  However, the applicability of this structure theorem is
bounded by the following facts:  Numerous module categories of interest fail
to be contravariantly finite;  on top of it, deciding whether a given
subcategory
$\A$ of $\lamod$ has this property is a difficult task in general.  And even
when the question of contravariant finiteness has been resolved in the
positive, it may still be extremely challenging to pin down the minimal right
$\A$-ap\-prox\-i\-ma\-tions of the simple objects.

In response to these obstacles, we suggest three lines of investigation, all
based on the concept of approximation or derivatives thereof.  The present
report is partly an overview over existing results in these directions,
illustrated by examples, partly an extension of ideas in the literature (e.g.,
most of Section 8 is new), and partly a program to be pursued more fully in
the future.  Roughly, the three lines of approach are as follows.

(1)  Find manageable positive and negative criteria for contravariant
finiteness, and systematically broaden the classes of algebras and module
categories for which the `classical' approximation theory can be brought
to bear.  In particular, this calls for an enlarged arsenal of
techniques to find and describe minimal
$\A$-ap\-prox\-i\-ma\-tions in the case  of contravariant finiteness and for a study of
`typical shapes' occurring among the minimal approximations of the simple
modules, in dependence  of $\la$ and $\A$.  This line is illustrated by
Sections 3, 4, and 7 below.

(2)  Given a class $\C$ of representations which cannot be mastered in the
classical mode of direct sum decompositions and a complete listing of the
indecomposable summands occurring, scan the objects of $\C$ by comparing
them with the objects of a reference class which is tailored to measure for
this purpose.  More precisely, construct contravariantly finite subcategories
$\A$ of $\lamod$ so that the following classification of $\C$ relative to
$\A$ yields useful information on $\C$:  namely, modules
$C$ and $D$ are considered `similar relative to $\A$' provided they
have isomorphic minimal right $\A$-ap\-prox\-i\-ma\-tions.  This, usually coarser,
alternative to the traditional classification in terms of precise
structural data naturally
yields information of varying degrees of precision, ranging from the very rough
-- if $\A$ is small -- to a complete picture -- if $\C \subseteq \A$.  

Here we
illustrate this idea for
$\C = \pinf(\lamod)$, the category of all modules of finite projective
dimension in $\lamod$.  The objects of the reference category we use are
glued together from uniserial building blocks (Sections 5,6,8).  Sections
5 and 6 provide the tools:  The first addresses the classification of the
uniserial left $\la$-modules in an informal sketch of some of the author's
results in \cite{Geom I,II,III}.  The second proposes a pattern for pasting
uniserial modules together, while keeping a tight grip on the
structure of the  new objects created; our choice of pattern is justified in
Section 7.  In Section 8, we deal with algebras of finite uniserial type or,
more precisely, with algebras for which the reference categories constructed in
Section 7 have finite representation type.  We conclude with
concrete instances of the outlined approximation strategy.

(3)  In a variation of our theme, we describe infinite dimensional
substitutes for minimal $\A$-ap\-prox\-i\-ma\-tions as developed by the
author and Happel in \cite{10}.  These encode the same type of information
as finite dimensional approximations do when they exist.  To outline the
underlying idea we require a few concepts, the first of which
is only a slight extension of the original definition of right
$\A$-ap\-prox\-i\-ma\-tions.  If $\B
\subseteq \A$ and $M$ is a left $\la$-module, a map $f: A \rightarrow M$ is
called a (right)
$\B$-{\it approximation of} $M$ {\it inside} $\A$ if $A \in \A$ and each
map $g \in$ Hom$_{\la}(B,M)$ for $B\in \B$ factors through $f$.  Furthermore, a
$\la$-module $H$ is called an $\A$-{\it phantom} of $M$ if there exists a
nonempty finite subclass $\B$ of $\A$ such that each
$\B$-ap\-prox\-i\-ma\-tion of $M$ inside $\A$ has $H$ as a subfactor;
direct limits of such modules
$H$ -- not required to be finitely generated -- are again named phantoms.  
We are particularly interested in phantoms $H_j$ with the property that all
homomorphisms $B \rightarrow M$ with $B$ running through a specified subclass
$\B_j$ of $\B$ factor through $H_j$.  Given a chain $\B_1 \subseteq \B_2
\subseteq \B_3 \subseteq \dots$ of subclasses of $\B$, the corresponding
sequence of phantoms
$H_1, H_2, H_3, \dots$ will not stabilize in general, and so a direct limit is
the natural way to integrate the information represented by such a
sequence.  To visualize this type of information, suppose that $M = S$ is
simple; then a phantom $H_j$ as above communicates -- in the most compressed
form possible within the category $\A$ -- the relations of those objects $B \in
\B_j$ which contain $S$ in their tops.  

The key fact in this connection is that nontrivial phantoms always exist.  In
fact, a module
$M \in \lamod$ fails to have an $\A$-ap\-prox\-i\-ma\-tion in the traditional
sense if and only if
$M$ has $\A$-phantoms of countably infinite
vector space dimension.

\medskip 

In Section 2 we review the basics of contravariantly finite
subcategories, concentrating on the point of view to be pursued later.  In
particular, we will omit all results linking functorially finite
categories to tilting theory or relative derived functors.  Since the
concept of covariant finiteness of a subcategory $\A$ of $\lamod$ is dual
to that of contravariant finiteness, we restrict our attention to the
latter; the lines of investigation we discuss can be dualized accordingly. 
The following convention will therefore not lead to ambiguities: namely, we
will refer to {\it right} $\A$-{\it approximations} simply as $\A$-{\it
approximations}. 

\medskip

Throughout, $\la$ will be a finite dimensional algebra over a field
$K$ with Jacobson radical $J$. We assume $\la$ to be split, i.e., $\la$ is of
the form
$\kgam / I$, where
$\Gamma$ is a quiver and $I$ is an admissible ideal of the path algebra
$\kgam$.  Our convention for composing paths in $\kgam$ is as follows:  if $p$
and $q$ are paths, then $pq$ stands for `$p$ after $q$' if the endpoint of $q$
coincides with the starting point of $p$, and $pq = 0$ otherwise.  A path
$u$ is called a {\it right subpath} (respectively, a {\it left subpath} of
$p$ in case there exists a path $v$ such that $p = vu$ (respectively, $p =
uv$); both left and right subpaths will also be referred to as {\it
subpaths} of $p$.
 	
Whenever we mention primitive idempotents of $\la$, we mean those associated
with the vertices of $\Gamma$;  in fact, we identify these vertices with
the corresponding idempotents.  Given a left $\la$-module $M$, a {\it top
element} of $M$ is an element $x \in M \setminus JM$ such that $x = ex$
for some primitive idempotent $e$;  in that case we also say that $x$ is
a {\it top element of type} $e$ of $M$.

Like numerous other authors (see, e.g., \cite{1} and \cite{7}), we find it
convenient to represent certain
$\la$-modules by graphs. A special breed
of such graphs, based on sequences of top elements of the pertinent modules,
will make a great deal of information available at a glance.  We introduce
these labeled and layered graphs informally by
means of some illustrative examples.

\example{Example 1.1} Let $\la=\kgam/I$, where $\gam$ is the quiver

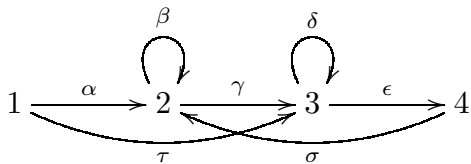
\ignore{
$$\xymatrixcolsep{3pc}
\xy\xymatrix{
1 \arrow[r]^\alpha \arrow@/_1pc/[rr]_\tau &2 \arrow[r]^\gamma \uloopr{\beta}
&3 \arrow[r]^\epsilon \uloopr{\delta} &4 \arrow@/^1pc/[ll]^\sigma
}\endxy$$
}

\noindent and $I$ an ideal in $\kgam$. That a left $\la$-module $M$ have
(layered and labeled) graph

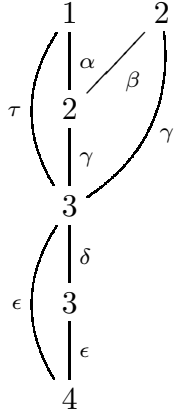
\ignore{
$$\xymatrixcolsep{1.5pc}\xymatrixrowsep{1.5pc}
\xy\xymatrix{
1 \edge[d]^\alpha \edge@/_1pc/[dd]_\tau &2 \edge[dl]^\beta
\edge@/^1pc/[ddl]^\gamma\\
2 \edge[d]^\gamma\\
3 \edge[d]^\delta \edge@/_1pc/[dd]_\epsilon\\
3 \edge[d]^\epsilon\\
4}\endxy$$
}

\noindent relative to top elements $x_1$ of type $e_1$ and $x_2$ of type
$e_2$ means that:
$$\gather M/JM\cong S_1\oplus S_2; \qquad JM/J^2M\cong S_2;\\
J^2M/J^3M\cong J^3M/J^4M\cong S_3; \qquad J^4M\cong S_4;\\
\alpha x_1 \text{\ and\ }
\beta x_2 \text{\ each generate\ } JM \text{\ modulo\ } J^2M;\\
\tau x_1,\ \gamma x_2,\ \gamma\alpha x_1 \text{\ (and hence also\ }
\gamma\beta x_2 \text{) each generate\ } J^2M \text{\ modulo\ } J^3M;\\
\delta\gamma\alpha x_1 \text{\ generates\ } J^3M \text{\ modulo\ } J^4M;\\
\epsilon\gamma\alpha x_1 \text{\ and\ } \epsilon\delta\gamma\alpha x_1
\text{\ each generate\ } J^4M. \endgather$$
Finally, whenever $p$ is a path in $\gam$ starting in $e_i$ which is not
recorded under $x_i$ in the graph of $M$, the product $px_i$ is zero.

On the other hand,

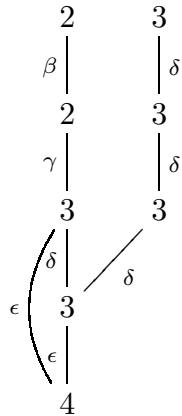
\ignore{
$$\xymatrixcolsep{1.5pc}\xymatrixrowsep{1.5pc}
\xy\xymatrix{
2 \edge[d]_\beta &3 \edge[d]^\delta\\
2 \edge[d]_\gamma &3 \edge[d]^\delta\\
3 \edge[d]_\delta \edge@/_1pc/[dd]_\epsilon &3 \edge[dl]^\delta\\
3 \edge[d]_\epsilon\\
4}\endxy$$
}

\noindent is not the graph of a left $\la$-module in case $\epsilon\delta^3=0$
in
$\la$.

It is often convenient to communicate the ideal $I$ by giving the graphs of
the indecomposable projective modules $\la e_i$. The algebra $\la$ with
quiver $\gam$ as above and with the following indecomposable projectives will
recur in Section 6.

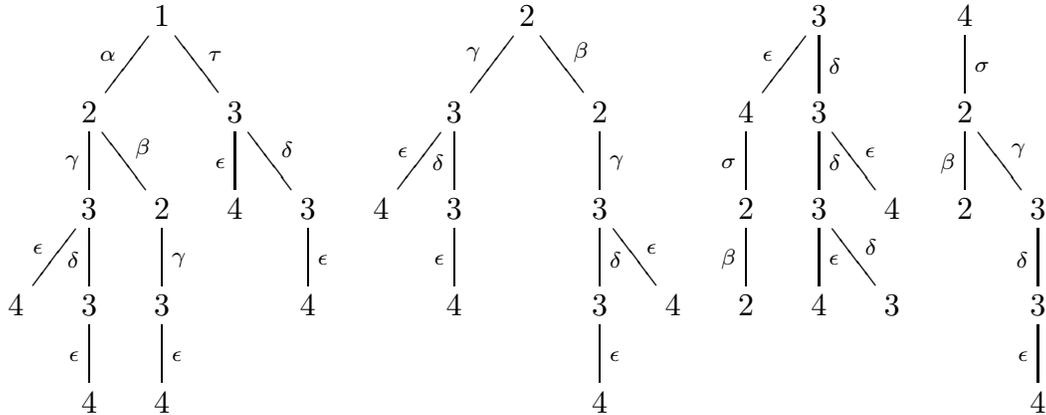
\ignore{
$$\xymatrixcolsep{1pc}\xymatrixrowsep{1.5pc}
\xy\xymatrix{
 &&1 \edge[dl]_\alpha \edge[dr]^\tau &&&&&2 \edge[dl]_\gamma \edge[dr]^\beta 
&&&&3 \edge[dl]_\epsilon \edge[d]^\delta &&4 \edge[d]^\sigma\\
 &2 \edge[d]_\gamma \edge[dr]^\beta &&3 \edge[d]_\epsilon \edge[dr]^\delta
&&&3 \edge[dl]_\epsilon \edge[d]_\delta &&2 \edge[d]^\gamma &&4
\edge[d]_\sigma &3 \edge[d]^\delta \edge[dr]^\epsilon &&2 \edge[d]_\beta
\edge[dr]^\gamma\\
 &3 \edge[dl]_\epsilon \edge[d]_\delta &2 \edge[d]^\gamma &4 &3
\edge[d]^\epsilon &4 &3 \edge[d]_\epsilon &&3
\edge[d]^\delta \edge[dr]^\epsilon &&2 \edge[d]_\beta &3 \edge[d]^\epsilon
\edge[dr]^\delta &4 &2 &3 \edge[d]_\delta\\
4 &3 \edge[d]_\epsilon &3 \edge[d]^\epsilon &&4 &&4 &&3
\edge[d]^\epsilon &4 &2 &4 &3 &&3 \edge[d]_\epsilon\\
 &4 &4 &&&& &&4 &&&&&&4
}\endxy$$
}

We extend these conventions as follows: That $N\in
\lamod$ has the following graph

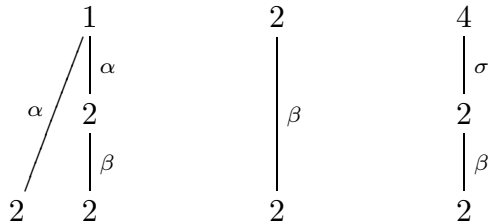
\ignore{
$$\xymatrixcolsep{1pc}\xymatrixrowsep{1.5pc}
\xy\xymatrix{
 &1 \edge[ddl]_\alpha \edge[d]^\alpha &&&2 \edge[dd]^\beta &&&4
\edge[d]^\sigma\\
 &2 \edge[d]^\beta &&& &&&2 \edge[d]^\beta\\
2 &2 &&&2 &&&2
}\endxy$$ 
}

\noindent is to signify that $J^2N=S_2^3$, and that there are top elements
$x_1$, $x_2$, $x_3$ of the obvious types such that $\alpha x_1$, $\beta\alpha
x_1$, $\beta x_2$, $\beta\sigma x_4$ generate $J^2N$, with every choice of
three of these four elements $K$-linearly independent (while the set of four
is $K$-linearly dependent). \qed\endexample

Finally, $\pinf(\lamod)$ will stand for the full subcategory of $\lamod$
having as objects the finitely generated modules of finite projective
dimension, and $\pinf(\Lambda\operatorname{-Mod})$ will be the analogous
subcategory of $\Lambda\operatorname{-Mod}$. As usual, $\lfd\la$ denotes the
supremum of the projective dimensions attained on $\pinf(\lamod)$.

\head 2. Contravariantly finite subcategories of $\lamod$\endhead

The notions of co- and contravariantly finite subcategories of $\lamod$ were
first introduced by Auslander and Smal\o\ in \cite{3}, in connection with
the existence of internal almost split sequences of a subcategory $\A$ and
preinjective/preprojective partitions of $\A$. The following definitions
evolved in \cite{3}, \cite{4} and \cite{2}. 

Throughout, we let $\A\subseteq \lamod$ be a full subcategory which is closed
under isomorphisms, direct summands, and finite direct sums.

\definition{Definitions 2.1} (1) Given $M\in\lamod$, a homomorphism $f :
A\rightarrow M$ with $A\in\A$ is called a {\it right
$\A$-ap\-prox\-i\-ma\-tion} of
$M$ if each $g\in \Hom_\la(B,M)$ with $B\in\A$ factors through $f$. Left
$\A$-ap\-prox\-i\-ma\-tions are defined dually.

(2) $\A$ is called {\it contravariantly finite} (resp., {\it covariantly
finite}) in $\lamod$ if each module $M\in\lamod$ has a right (resp., left)
$\A$-ap\-prox\-i\-ma\-tion.\enddefinition

In the following we will focus on contravariant finiteness and will briefly
write `$\A$-ap\-prox\-i\-ma\-tion' for `right $\A$-ap\-prox\-i\-ma\-tion'. In
\cite{2} it is shown that, whenever $M\in\lamod$ has an
$\A$-ap\-prox\-i\-ma\-tion, any two
$\A$-ap\-prox\-i\-ma\-tions of minimal $K$-dimension are isomorphic. We therefore
refer to {\it the minimal $\A$-ap\-prox\-i\-ma\-tion} of $M$ in that case.

The nomenclature in Definition 2.1(2) stems from the fact that $\A$ is
contravariantly finite in $\lamod$ if and only if the restricted
contravariant $\Hom$-functor $\Hom(-,M)|_\A : \A\rightarrow \bold{Ab}$ is
finitely generated for each object $M\in\lamod$. The flavor of this
finiteness condition of $\A$ relative to $\lamod$ is caught in the following
easy argument showing that each subcategory $\A\subseteq \lamod$ of finite
representation type is contravariantly (as well as covariantly) finite in
$\lamod$ \cite{AS, Proposition 4.2}: Indeed, if $A_1,\dots,A_m\in \A$ are
such that each object in $\A$ is a finite direct sum of copies of the $A_i$,
then given $M\in\lamod$, the homomorphism $\oplus_{i=1}^m \oplus_{j=1}^{d_i}
f_{ij} : \oplus_{i=1}^m A_i^{d_i} \rightarrow M$ is an
$\A$-ap\-prox\-i\-ma\-tion of
$M$ provided that, for each $i\le m$, the family $(f_{ij})_{j\le d_i}$ is a
$K$-basis of $\Hom_\la(A_i,M)$.

One of the first results to draw attention to these finiteness conditions for
subcategories of $\lamod$ was the following theorem of Auslander and Smal\o.

\proclaim{Theorem 2.2} \cite{4, Theorem 2.4} If $\A$ is covariantly and
contravariantly finite in $\lamod$, then $\A$ has almost split sequences.
\qed\endproclaim

The following two results, due to Auslander and Reiten, are pivotal in our
present investigation.

\proclaim{Theorem 2.3} \cite{2, Proposition 3.7} If $\A\subseteq \lamod$ is
a resolving subcategory (meaning that $\A$ contains the projectives in
$\lamod$ and is closed under extensions and kernels of epimorphisms), then
$\A$ is contravariantly finite in $\lamod$ if and only if each of the simple
left $\la$-modules has an $\A$-ap\-prox\-i\-ma\-tion. \qed\endproclaim

\proclaim{Theorem 2.4} \cite{2, Proposition 3.8} Suppose that $\A$ is a
resolving, contravariantly finite subcategory of $\lamod$. Moreover, suppose
that $\lamod$ contains $n$ isomorphism classes of simple modules and that
$A_1,\dots,A_n$ are the minimal $\A$-ap\-prox\-i\-ma\-tions of representative
simples. Then a left $\la$-module $X$ belongs to $A$ precisely when $X$ is a
direct summand of a module having a filtration with consecutive factors in
$\{A_1,\dots,A_n\}$. \qed\endproclaim

In the next section, we give a first instalment of `positive' examples, i.e.,
of instances of contravariantly finite subcategories. A powerful tool for
extending the list of such subcategories is the following theorem of Sikko
and Smal\o, which subsumes earlier results by Grecht, Vossieck, de la
Pena-Simson, Ringel, and Smal\o.

\proclaim{Theorem 2.5} \cite{22, Theorem 2.6} If $\A$ and $\B$ are
contravariantly finite subcategories of $\lamod$, then the category of all
direct summands of extensions of objects in $\A$ by objects in $\B$ has the
same property. \qed\endproclaim

\head 3. First instalment of positive examples\endhead

Roughly speaking, a subcategory $\A\subseteq \lamod$ is likely to be
contravariantly finite if it is either very small or very large. In the
former situation one usually obtains very rough approximations, in the latter
one can approach the target objects very closely with objects from $\A$. The
extreme cases are that of a category of finite representation type on one end
-- we already know such a category to be contravariantly finite in $\lamod$
-- and the case $\A=\lamod$ on the other end. While the latter does not hold
much interest, the former does. For example, in case the objects of $\A$ are
precisely the finitely generated projective left $\la$-modules, the minimal
$\A$-ap\-prox\-i\-ma\-tions coincide with projective covers. So in
particular, if
$\lfd\la=0$, the category $\pinf(\lamod)$ is contravariantly finite in
$\lamod$. The most interesting cases of contravariant finiteness of
$\pinf(\lamod)$, however, are those where
$\lfd\la>0$, while $\gldim\la =\infty$.

As was already observed by Auslander and Reiten in \cite{2}, $\pinf(\lamod)$
is contravariantly finite in $\lamod$ whenever $\la$ is stably equivalent to
a hereditary algebra. This includes the case where $J^2=0$, a situation in
which it is easy to describe the minimal $\pinf(\lamod)$-ap\-prox\-i\-ma\-tions.
According to \cite{12}, they look as follows: Given any $X\in\lamod$
with $JX=0$, write $X= X_{fin}\oplus X_{inf}$, where $X_{fin}$ is the sum of
those simple submodules of $X$ which have finite projective dimension, and
$X_{inf}$ is the sum of those simples which have infinite projective
dimension. Then, given any $M\in\lamod$ with projective cover $P\rightarrow
M$ say, the induced epimorphism $P/\Omega^1(M)_{fin} \rightarrow M$ is the
minimal $\pinf(\lamod)$-ap\-prox\-i\-ma\-tion of $M$. In particular, the
minimal
$\pinf(\lamod)$-ap\-prox\-i\-ma\-tion of a simple module $S= \la e/Je$ has
the form
$\la e/(Je)_{fin} \rightarrow S$. To be more specific, if $\la=
\kgam/\langle \text{all paths of length 2}\rangle$, where $\gam$ is the
quiver

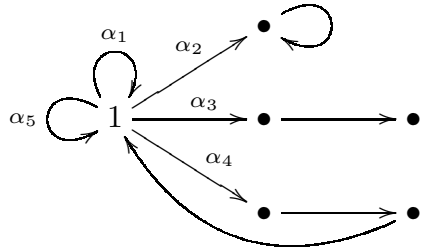
\ignore{
$$\xymatrixcolsep{3pc}\xymatrixrowsep{1.5pc}
\xy\xymatrix{
 &\bullet \rloopd{}\\
1 \arrow[ur]^(0.6){\alpha_2} \arrow[r]^(0.6){\alpha_3}
\arrow[dr]^(0.6){\alpha_4}
\uloopr{\alpha_1} \lloopd{\alpha_5} &\bullet \arrow[r] &\bullet\\
 &\bullet \arrow[r] &\bullet \arrow@/^2pc/[ull]
}\endxy$$
}

\noindent the graph of the minimal $\pinf(\lamod)$-ap\-prox\-i\-ma\-tion of
$S_1$ is the following `brush':

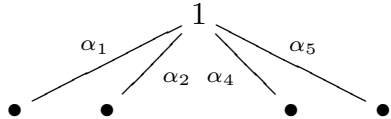
\ignore{
$$\xymatrixcolsep{1.5pc}\xymatrixrowsep{1.5pc}
\xy\xymatrix{
 &&1 \edge[dll]_{\alpha_1} \edge[dl]^{\alpha_2}
\edge[dr]_{\alpha_4} \edge[drr]^{\alpha_5}\\
\bullet &\bullet &&\bullet &\bullet
}\endxy$$
}

Moreover, Auslander and Reiten proved that for each Gorenstein algebra
$\la$, i.e., for each algebra $\la$ which has finite injective dimension on
both sides, the category $\pinf(\lamod)$ is contravariantly finite in
$\lamod$ (see \cite{2, p\. 150}). This result applies to a particularly
interesting situation to which the author was alerted by Kirkman and
Kuzmanovich: Namely, let $\O\subseteq M_n(F)$ be a tiled classical order
over a DVR $D$ with quotient field $F$; that $\O$ be tiled means that $\O$
contains a full set of $n$ primitive orthogonal idempotents. Moreover, let
$\pi$ be a uniformizing parameter of $D$ and observe that then $\la= \opio$
is a finite dimensional algebra over the residue class field $K=D/(\pi)$ of
$D$. If $\O$ has finite left and right injective dimension -- a fortiori, if
$\gldim \O<\infty$ -- one of the classical change of rings
theorems (see \cite{20, Theorem 205}) yields that $\la$ is a Gorenstein
algebra. So in particular $\pinf(\lamod)$ is contravariantly finite in
$\lamod$ whenever $\gldim\O<\infty$; hence, $\lfd\la=\sup \{\pdim A_i \mid
1\le i\le n\}$, where the $A_i$ are the minimal
$\pinf(\lamod)$-ap\-prox\-i\-ma\-tions of the simple left $\la$-modules. On the
other hand, $\lfd\la =\lfd\O-1 =\gldim\O-1$ if $\gldim\O<\infty$, due to
\cite{9}. Via approximation theory, this equality may provide access to
the possible values of $\gldim\O$, which have been the object of
investigation for a long time (see \cite{23}, \cite{18}, \cite{19}, \cite{21},
\cite{6} for more detail). A summary of the results to date pertaining to
the global dimensions of tiled classical orders over DVRs can be found in
\cite{13}.

\example{Problem 3.1} Let $\la=\opio$ for a tiled classical order $\O$ over a
DVR. Give an explicit description of the minimal
$\pinf(\lamod)$-ap\-prox\-i\-ma\-tions of the simple left $\la$-modules in terms of
the valuated quiver of $\O$ (see \cite{24}) or, equivalently, in terms of
quiver and relations of $\la$. \endexample

\head 4. Negative examples and a criterion for failure of contravariant
finiteness\endhead

The first example of a finite dimensional algebra $\la$ for which
$\pinf(\lamod)$ fails to be contravariantly finite is due to Igusa, Smal\o,
Todorov \cite{17}. It is a monomial relation algebra with vanishing radical
cube which, in addition, is biserial. (In Section 7, we will see that, by
contrast, $\pinf(\lamod)$ is always contravariantly finite in $\lamod$ when
$\la$ is left serial.) However, the conclusion that contravariant finiteness
of $\pinf(\lamod)$ in $\lamod$ be a rare occurrence would be precipitous. In
fact, the condition cuts diagonally across the standard groupings of finite
dimensional algebras, being extremely sensitive to changes in quiver and
relations. 

\proclaim{Criterion 4.1} \cite{10} Let $\A$ be a full subcategory of
$\lamod$. Moreover, suppose that $e_1,\dots,e_m$ are pairwise orthogonal
primitive idempotents of $\la$, and $p_1,\dots,\allowmathbreak
p_m,\allowmathbreak q_1,\dots,q_m$ elements of
$J$, with $p_i= p_ie_i$ and $q_i =q_ie_i$, such that the following conditions
are satisfied:

{\rm (1)} For each $n\in\N$, there is a module $M_n\in \A$, together with a
sequence $x_{n1},\dots,x_{nn}$ of top elements of $M_n$ which are
$K$-linearly independent modulo $JM_n$ and have the property that
$p_{r(i)}x_{ni}= q_{r(i+1)}x_{n,i+1}$ is nonzero for $1\le i\le n$, where
$r(i)\in
\{1,\dots,m\}$ is congruent to $i$ modulo $m$.

{\rm (2)} Whenever $A\in\A$, the following hold:

\quad{\rm (i)} If $x\in A$ is a top element of type $e_1$, then $p_1x\ne0$,
and

\quad{\rm (ii)} If $y,z\in A$ are such that $p_{r(i)}y= q_{r(i+1)}z \ne 0$,
then
$p_{r(i+1)}z \ne 0$.

\noindent It follows that $S_1= \la e_1/Je_1$ does not have an
$\A$-ap\-prox\-i\-ma\-tion.
\qed\endproclaim

\example{Example 4.2} \cite{17} Let $\A=\pinf(\lamod)$, where $\la=\kgam/I$
is the monomial relation algebra with quiver 

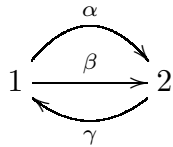
\ignore{
$$\xymatrixcolsep{3pc}
\xy\xymatrix{
1 \arrow[r]^\beta \arrow@/^1.5pc/[r]^\alpha &2 \arrow@/^1pc/[l]^\gamma
}\endxy$$
}

\noindent such that the indecomposable projective left $\la$-modules have
graphs 

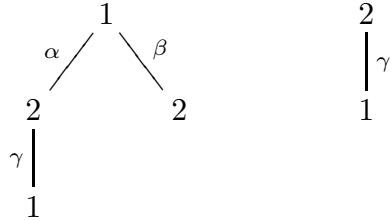
\ignore{
$$\xymatrixcolsep{1pc}\xymatrixrowsep{1.5pc}
\xy\xymatrix{
 &1 \edge[dl]_\alpha \edge[dr]^\beta &&&&2 \edge[d]^\gamma\\
2 \edge[d]_\gamma &&2 &&&1\\
1}\endxy$$
}

Set $m=1$, $p_1=\beta$, $q_1=\alpha$, and $S_i=\la e_i/Je_i$. Then the
modules $M_n\in \pinf(\lamod)$ with graphs

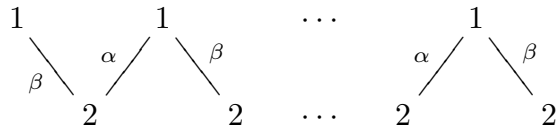
\ignore{
$$\xymatrixcolsep{1pc}\xymatrixrowsep{1.5pc}
\xy\xymatrix{
1 \edge[dr]_\beta &&1 \edge[dl]_\alpha \edge[dr]^\beta &&\cdots &&1
\edge[dl]_\alpha \edge[dr]^\beta\\
 &2 &&2 &\cdots &2 &&2
}\endxy$$
}

\noindent having $n$ linearly independent top elements of type $e_1$ modulo
the radical clearly satisfy Condition (1) of Criterion 4.1.

Moreover, Condition (2i) of the criterion is trivially satisfied, in view of
the fact that
$\pdim S_2 =\infty$; indeed, given any $A\in \lamod$ with top element
$x=e_1x$ and $p_1x= \beta x=0$, the simple module $S_2$ is a direct summand
of $\Omega^1(A)$. As for (2ii), whenever $p_1y= \beta y =q_1z =\alpha z
\ne0$, the element $e_1z$ is a top element of $A$ of type $e_1$. Thus, $S_1$
does not have a $\pinf(\lamod)$-ap\-prox\-i\-ma\-tion. \qed\endexample

To the following example we will refer back in Sections 8 and 9. It is an
instance where Criterion 4.1 applies only for $m>1$.  

\example{Example 4.3} This time, let $\la=\kgam/I$, where $\gam$ is the quiver

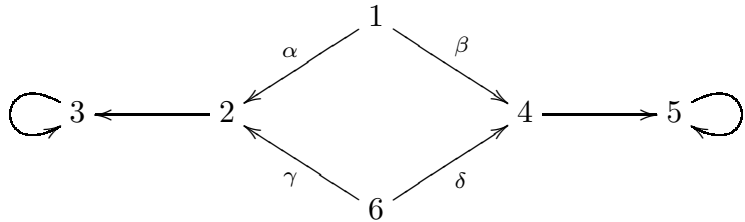
\ignore{
$$\xymatrixcolsep{3pc}\xymatrixrowsep{1.5pc}
\xy\xymatrix{
 &&1 \arrow[dl]_\alpha \arrow[dr]^\beta\\
3 \lloopd{} &2 \arrow[l] &&4 \arrow[r] &5 \rloopd{}\\
 &&6 \arrow[ul]^\gamma \arrow[ur]_\delta
}\endxy$$
}

\noindent and $I\subset\kgam$ is such that the indecomposable projective left
$\la$-modules have the following graphs:

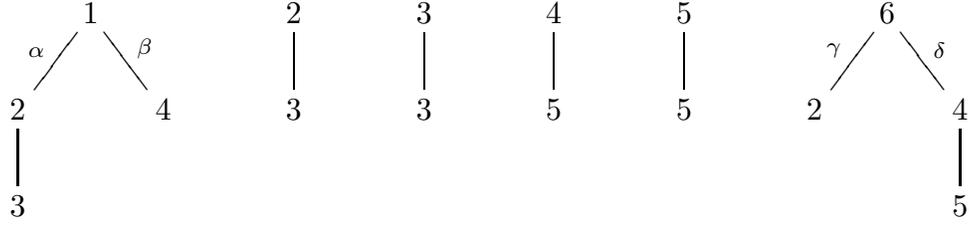
\ignore{
$$\xymatrixcolsep{1pc}\xymatrixrowsep{1.5pc}
\xy\xymatrix{
&1 \edge[dl]_\alpha \edge[dr]^\beta &&&2 \edge[d] &&3 \edge[d] &&4 \edge[d] &&5
\edge[d] &&&6 \edge[dl]_\gamma \edge[dr]^\delta\\
2 \edge[d] &&4 &&3 &&3 &&5 &&5 &&2 &&4 \edge[d]\\
3 &&&&& &&&&& &&&&5
}\endxy$$
}

Set $m=2$, $p_1=\beta$, $p_2=\gamma$, $q_1=\alpha$, $q_2=\delta$, and for
$n\in\N$, let $M_n\in\pinf(\lamod)$ be the module with graph

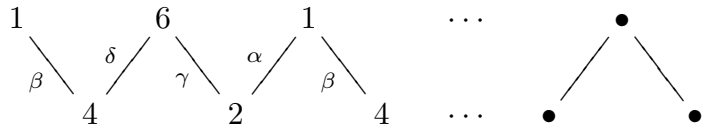
\ignore{
$$\xymatrixcolsep{1pc}\xymatrixrowsep{1.5pc}
\xy\xymatrix{
1 \edge[dr]_\beta &&6 \edge[dl]_\delta \edge[dr]_\gamma &&1 \edge[dl]_\alpha
\edge[dr]_\beta &&\cdots &&\bullet \edge[dl] \edge[dr]\\
 &4 &&2 &&4 &\cdots &\bullet &&\bullet
}\endxy$$
}

\noindent having $n$ linearly independent top elements modulo the radical, of
types alternating between 1 and 6. It is straightforward to see that Condition
(2) of Criterion 4.1 is satisfied as well, so that $\pinf(\lamod)$ again
fails to be contravariantly finite in $\lamod$. \qed\endexample

Our two final examples illustrate the instability -- under minor
modifications of quiver and/or relations -- of the condition that
$\pinf(\lamod)$ be contravariantly finite in $\lamod$. The first is a variant
of the Igusa-Smal\o-Todorov example (4.2); the second results from a further
slight alteration of the relations.

\example{Example 4.4} Let $\la=\kgam/I$ where $\gam$ is 

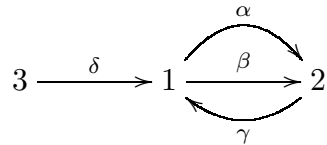
\ignore{
$$\xymatrixcolsep{3pc}
\xy\xymatrix{
3 \arrow[r]^\delta &1 \arrow[r]^\beta \arrow@/^1.5pc/[r]^\alpha &2
\arrow@/^1pc/[l]^\gamma
}\endxy$$
}

\noindent and $I$ is such that the indecomposable projective left
$\la$-modules have graphs

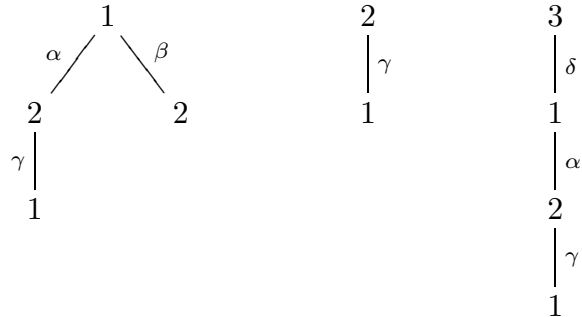
\ignore{
$$\xymatrixcolsep{1pc}\xymatrixrowsep{1.5pc}
\xy\xymatrix{
 &1 \edge[dl]_\alpha \edge[dr]^\beta &&&&2 \edge[d]^\gamma &&&3
\edge[d]^\delta\\ 
2 \edge[d]_\gamma &&2 &&&1 &&&1 \edge[d]^\alpha\\
1 && &&& &&&2 \edge[d]^\gamma\\
 && &&& &&&1
}\endxy$$
}

\noindent This algebra $\la$ contains the algebra of Example 4.2, and the
modules $M_n$ defined in 4.2 remain modules of finite projective dimension
over the new algebra. However, this time $\pinf(\lamod)$ is contravariantly
finite in $\lamod$; see \cite{10}. \qed\endexample

\example{Example 4.5} The quiver of the algebra $\la$ is the same as that of
Example 4.4, but we delete one of the relations in the previous example, to
the effect that the indecomposable projective left $\la$-modules take on the
forms

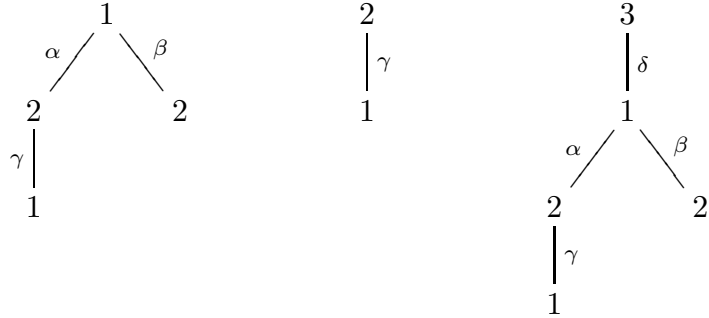
\ignore{
$$\xymatrixcolsep{1pc}\xymatrixrowsep{1.5pc}
\xy\xymatrix{
 &1 \edge[dl]_\alpha \edge[dr]^\beta &&&&2 \edge[d]^\gamma &&&&3
\edge[d]^\delta\\ 
2 \edge[d]_\gamma &&2 &&&1 &&&&1 \edge[dl]_\alpha \edge[dr]^\beta\\
1 && &&& &&&2 \edge[d]^\gamma &&2\\
 && &&& &&&1
}\endxy$$
}

\noindent Again, Criterion 4.1 (with $m=1$, $p_1=\beta$ and $q_1=\alpha$)
readily yields that
$S_1$ fails to have a
$\pinf(\lamod)$-ap\-prox\-i\-ma\-tion. \qed\endexample

\head 5. First intermezzo: Uniserial representations\endhead

We give a rough sketch of results from \cite{14,15,16} which will
provide the foundation for the construction of several useful contravariantly
finite subcategories of $\lamod$, the objects of which are well understood.
The pivotal problems addressed are the following:

(I) Classify the uniserial left $\la$-modules in terms of manageable
isomorphism invariants.

(II) Characterize the split algebras of finite uniserial type, i.e.,
characterize those algebras $\la=\kgam/I$ for which there are only
finitely many uniserial
$\la$-modules up to isomorphism.

As the author learned in the meantime, Auslander has proposed these problems
since the mid-70's. We include an excerpt of an email message of July 1993
from him to the author. ``I have been raising the question of the
classification or description of uniserial modules for artin algebras for
many years now. I believe the first time I raised the question in public was
at the special session on representation theory at the winter meeting of the
AMS in Atlanta around 1975. The big shot group representation people, like
[\dots] assured me that they would have an answer by the afternoon. I am
still waiting. I have raised the question repeatedly since then [\dots] As to
motivation. One reason I am interested in these modules, aside from the fact
that they should in some sense be the simplest nonprojective modules, is the
fact that they have bounded length for a given algebra. Therefore the second
Brauer-Thrall is true for algebras with an infinite number of uniserial
modules. Hence it would be interesting to know for which algebras there are
only a finite number [\dots] Secondly, I am interested where the uniserial
modules occur in AR-quivers and preprojective partitions since as you
mentioned other modules can be in some sense approximated by uniserial
modules. [\dots]''

The placement of uniserial modules in the Auslander-Reiten quiver will not be
discussed in this intermezzo. We just mention that Axel Boldt is working on
this subject in his dissertation, and has already settled the question in the
hereditary case.

A preliminary subdivision of the class of uniserial left $\la$-modules is
in terms of their `masts'.

\definition{Definition 5.1} Given a uniserial left $\la$-module $U$ of length
$l+1$, any path $p$ of length $l$ in $\kgam$ with $pU\ne0$ is called a {\it
mast} of $U$. \enddefinition

Observe that, if $U$ is uniserial with consecutive simple composition factors
$J^iU/J^{i+1}U\cong \la e(i)/Je(i)$ for $0\le i\le l$, and with layered and
labeled graph $G$, then the masts of $U$ correspond precisely to those edge
paths in $G$ which pass exactly once through each of these vertices from top to
bottom. E.g., a uniserial module with graph 

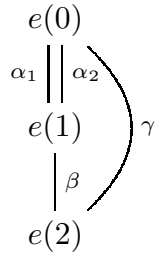
\ignore{
$$\xymatrixcolsep{1.5pc}\xymatrixrowsep{1.5pc}
\xy\xymatrix{
e(0) \edge[d]<-0.5ex>_{\alpha_1} \edge[d]<0.5ex>^{\alpha_2}
\edge@/^2pc/[dd]^\gamma\\
e(1) \edge[d]^\beta\\
e(2)}\endxy$$
}

\noindent has masts $\beta\alpha_1$ and $\beta\alpha_2$.

The starting point of our approach to uniserial representations is the
following theorem which we state somewhat informally.

\proclaim{Theorem 5.2} \cite{14,15} Let $p\in\kgam$ be a path.

{\rm (1)} There is an affine algebraic variety $V_p$, not necessarily
irreducible, which parametrizes the isomorphism types of the uniserial left
$\la$-modules with mast $p$ in a natural fashion. Somewhat more precisely,
there exists a canonical surjection
$$\Phi_p : V_p \rightarrow \{\text{isomorphism types of uniserials in
$\lamod$ with mast\ } p\}.$$
In particular, $V_p$ is nonempty if and only if there exists a uniserial left
$\la$-module with mast $p$. (Polynomials for this variety can be readily
determined on the basis of a `coordinatization' $\la=\kgam/I$.)

{\rm (2)} If $p : e(0) \rightarrow e(1) \rightarrow \cdots \rightarrow e(l)$
does not have a right subpath of positive length from $e(0)$ to $e(0)$, the
map $\Phi_p$ is bijective, and the points of $V_p$ serve as isomorphism
invariants of the uniserial modules with mast $p$. More sharply, if $e(0)$
recurs $t$ times among the vertices $e(1),\dots,e(l)$, then each fibre of
$\Phi_p$ is contained in a closed subvariety of $V_p$ of dimension at most
$t$.

{\rm (3)} If $V_p\subseteq \Bbb A^d$ and $t$ is as under {\rm (2)}, there
exists a system of equations
$$S_p(X,Y,Z) = S_p(X_1,\dots,X_d,Y_1,\dots,Y_d,Z_1,\dots,Z_t),$$
linear in $Z_1,\dots,Z_t$, with the following property: Two points $k,k'\in
V_p$ belong to the same fibre of $\Phi_p$ if and only if the linear system
$S_p(k,k',Z)$ is consistent. (Just as $V_p$, the system $S_p(X,Y,Z)$ can be
effectively computed from $\gam$ and $I$.)\qed\endproclaim

The traditional varieties of $\la$-modules of a fixed
$K$-dimension (or of the cyclics in $\lamod$) clearly contain the collection
of uniserial modules as open subvarieties. However, these varieties are far
too large and unwieldy to serve the purpose of an effective classification.
By contrast, the varieties which we will consider here fit the collection of
uniserials with a fixed composition series rather tightly. In fact, the bit
of slack which may occur in the presence of certain oriented cycles is quite
harmless; each such cycle just adds a copy of $\Bbb A^1$ to the pertinent
variety.

We will not give the equations defining the varieties $V_p$ here, but will
instead give a very rough answer to the question of what data are recorded by
the coordinates of their points.
Suppose that
$p=\alpha_l\cdots
\alpha_1$, where the
$\alpha_i$ are arrows. Given any uniserial $U\in\lamod$ with mast $p$ and a
top element
$x$, the products $\alpha_m\cdots \alpha_1x$, $0\le m\le l$, clearly form a
$K$-basis of $U$. Roughly speaking, the points of $V_p$ corresponding to $U$
are strings of coordinate vectors of elements $qx$, relative to this basis,
where
$q$ runs through certain paths in $\kgam$ (sufficiently many to pin down $U$
up to isomorphism). In particular, the graphs of uniserials are available at
a glance from the corresponding points on the varieties $V_p$. For details, we
refer to
\cite{14,15}. We include a fairly transparent example to illustrate the
correspondence of Theorem 5.2, again suppressing all technical detail.

\example{Example 5.3} Let $\gam$ be the
quiver

\ignore{
$$\xymatrixcolsep{1.5pc}
\xy\xymatrix{
1 \arrow[r]_{\alpha_1} \arrow@/^1pc/[rr]^{\gamma_1} &2
\arrow[r]_{\beta_1} &3
\arrow[r]_{\alpha_2} \arrow@/^1pc/[rr]^{\gamma_2} &4 
\arrow[r]_{\beta_2} &5
\arrow[r]_{\alpha_3} \arrow@/^1pc/[rr]^{\gamma_3} &6
\arrow[r]_{\beta_3} &7 
\arrow[r]_{\alpha_4} \arrow@/^1pc/[rr]^{\gamma_4} &8
\arrow[r]_{\beta_4} &9
\arrow[r]_{\alpha_5} \arrow@/^1pc/[rr]^{\gamma_5} &10 
\arrow[r]_{\beta_5} &11
}\endxy$$
}

\noindent and $p$ the path $q_5q_4q_3q_2q_1$, where $q_i=\beta_i\alpha_i$ for
$1\le i\le 5$. Define $\la =K\Gamma/I$, where $I$ is the ideal generated
by the following relations:
$$\gather \gamma_5\gamma_4q_3q_2q_1 -q_5q_4\gamma_3\gamma_2\gamma_1
+q_5q_4q_3q_2\gamma_1,\hskip0.4truein q_5q_4q_3q_2\gamma_1
-q_5q_4q_3\gamma_2q_1,\\
 q_5q_4q_3q_2\gamma_1 -q_5q_4\gamma_3q_2q_1,\hskip0.4truein
q_5\gamma_4q_3q_2q_1 -\gamma_5q_4q_3q_2q_1.\endgather$$ 
Then $V_p =V(X_{21}X_{22} -X_{11}X_{12}X_{13}+X_{11},\ X_{11}-X_{12},\
X_{11}-X_{13},\ X_{21}-X_{22})$. Observe that $V\cong V(X_2^2- X_1(X_1^2-1))$
is the elliptic curve with $\Bbb R$-graph

\ignore{
$$\xy (0,0);(0,0)
\curve{(0,2)&(-2,7)&(-10,7)&(-10,-7)&(-2,-7)&(0,-2)}
\endxy
\hskip1.0truecm\xy (35,15);(35,-15)
\curve{(30,10)&(18,7)&(10,4)&(10,-4)&(18,-7)&(30,-10)}
\endxy$$
}
\endexample

A host of questions supplementing Problems I and II pose themselves at this
point. We will briefly address a few of the most immediate among these.

\example{Additional Questions 5.4} {\rm (1)} Suppose that $p,q\in \kgam$ are
two paths of length $l$, both passing through the sequence of vertices
$(e(0),\dots,e(l))$. Then the uniserials with mast $p$ and those with mast
$q$ have the same sequence of composition factors, namely $(\la
e(0)/Je(0),\dots, \la e(l)/Je(l))$. What can be said about the intersection
$\Phi_p(V_p)\cap \Phi_q(V_q)$?

{\rm (2)} The varieties $V_p$ are defined by polynomials over $K$ which
depend on $\gam$ and $I$; the labeling as such is, in fact, tied to the
coordinatization of $\la$. Is there geometric information on the uniserial
left $\la$-modules with fixed sequence of composition factors which does not
depend on the coordinatization of $\la$?

{\rm (3)} Which affine varieties arise as varieties of uniserial
modules with a fixed sequence of composition factors?\endexample

\example{Answers (sketch)} {\rm (1)} \cite{14, Theorem D} Given any pair
of irreducible components $U_p\subseteq V_p$ and $U_q\subseteq V_q$, the
intersection $D= \Phi_p(U_p)\cap \Phi_q(U_q)$ is either empty or else
$\Phi_p^{-1}(D)$ and $\Phi_q^{-1}(D)$ are dense open subsets of $U_p$ and
$U_q$ respectively, and there exists an isomorphism of varieties $\Psi :
\Phi_p^{-1}(D) \rightarrow \Phi_q^{-1}(D)$ which makes the following diagram
commutative:

\ignore{
$$\xy\xymatrix{ \Phi_p^{-1}(D) \arrow[rr]^\Psi \arrow[dr]_{\Phi_p}
&&\Phi_q^{-1}(D) \arrow[dl]^{\Phi_q}\\
 &D}\endxy$$
}

\noindent In particular, $U_p$ and $U_q$ are birationally equivalent in the
latter case. 

{\rm (2)} \cite{16} Let $\Bbb S= (S(0),\dots,S(l))$ be a sequence of simple
left
$\la$-modules.

If $\gam$ has no double arrows, there is at most one path $p$ of length $l$
passing through a sequence of primitive idempotents of $\la$ corresponding to
the simples $S(i)$ and, in the positive case,
the variety
$V_p =V_{\Bbb S}$ is determined up to isomorphism by the
$K$-algebra isomorphism type of $\la$.

In the general case, start with a coordinatization $\la=\kgam/I$ and let
$(e(0),\dots,\allowmathbreak e(l))$ be the sequence of vertices of $\gam$
corresponding to the simple modules $S(i)$. Denote by $V_{\Bbb S}$ the set of
the birational equivalence classes of the irreducible components of the
$V_p$'s, where $p$ runs through the paths of length $l$ passing through
$(e(0),\dots,e(l))$. Then $V_{\Bbb S}$ is uniquely determined by the
isomorphism type of $\la$.

{\rm (3)} Each affine variety occurs as a variety $V_p =V_{\Bbb S}$, even
under the additional requirements that the map $\Phi_p$ be bijective and that
$V_p$ be determined up to isomorphism by the corresponding algebra $\la$. More
precisely: Given any affine algebraic variety $V$
over $K$, there exists an acyclic quiver $\Gamma$ without double arrows,
together with a path $p$ in $K\Gamma$, such that $V\cong V_p$. \qed\endexample

To close in on the structure of the algebras of finite uniserial type, we
require several additions to our conceptual framework (see \cite{15}).
Since, in the sequel, we will focus on algebras with the stronger property
that all the varieties $V_p$ be finite, we will content ourselves with the far
more straightforward characterization of these latter algebras. This
characterization will be preceded by a strong necessary condition for finite
uniserial type.

\proclaim{Theorem 5.5} \cite{15} If $\la$ has finite uniserial type, then
the following condition (N) is satisfied: Whenever
$\alpha : e\rightarrow e'$ is an arrow in $\Gamma$ and $p : e\rightarrow e'$
a mast of positive length, the path $p$ belongs to $K\Gamma\alpha \cup \alpha
K\Gamma$, that is,
$p$ is of the form

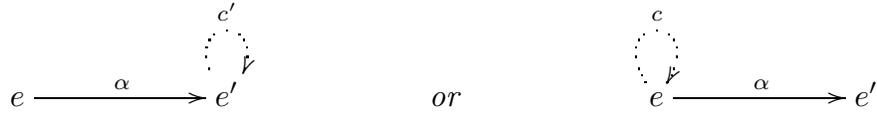
\ignore{
$$\xy\xymatrix{
e \arrow[rr]^\alpha &&e' \udotloopr{c'} &&\text{or}&& e \arrow[rr]^\alpha
\udotloopr{c} &&e'
}\endxy$$
}

\noindent where $c',c$ are oriented cycles which may be trivial.

A coordinate-free rendering of condition (N) is as follows: If there exists
a uniserial left $\la$-module $W$ of length 2 with top $S$ and socle $S'$,
and if $U$ is any uniserial left $\la$-module of length $l\ge2$ with top $S$
and socle
$S'$, then either $U/J^2U\cong W$ or else $J^{l-2}U\cong W$.

Condition (N), in turn, implies that:

{\rm (a)} $\Gamma$ has no double arrows, meaning that each uniserial left
$\la$-module has a unique mast.

{\rm (b)} Given any uniserial module $U$ with mast $p$, each graph of $U$
results from the superposition of graphs of the form

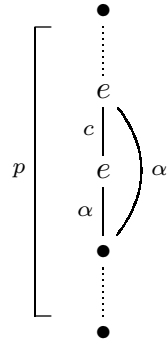
\ignore{
$$\xymatrixrowsep{1.25pc}
\xy\xymatrix{
\bullet \arrow@{.}[d] \arrow@^{|-|}[dddd]<-5ex>_p\\
e \edge[d]_c \edge@/^1pc/[dd]^\alpha\\
e \edge[d]_\alpha\\
\bullet \arrow@{.}[d]\\
\bullet}\endxy$$
}

\noindent under identification of the edge path $p$; here $c$ is an oriented
cycle of positive length and $\alpha$ an arrow such that $\alpha c$ is a
subpath of $p$. \qed\endproclaim

The algebras we wish to describe turn out to be precisely those for which the
second option in Condition (N) of Theorem 5.5 is excluded.

\proclaim{Theorem 5.6} \cite{15} Given any algebra $\la= \kgam/I$, the
following statements are equivalent:

{\rm (1)} For each path $p\in K\Gamma$,  the variety $V_p$ is finite.

{\rm (2)} For each path $p\in K\Gamma$, the variety $V_p$ is either empty or
a singleton.

{\rm (3)} Whenever $\alpha : e\rightarrow e'$ is an arrow in $\Gamma$ and $p
: e\rightarrow e'$ a mast of positive length, the path $p$ is equal
to $c'\alpha$, where
$c'$ is an oriented cycle which may be trivial, i.e., $p$ has the form

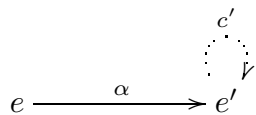
\ignore{
$$\xy\xymatrix{
e \arrow[rr]^\alpha &&e' \udotloopr{c'} 
}\endxy$$
}

In coordinate-free terms: If there exists
a uniserial left $\la$-module $W$ of length 2 with top $S$ and socle $S'$,
and if $U$ is any uniserial left $\la$-module of length $\ge2$ with top $S$
and socle
$S'$, then $U/J^2U\cong W$.

{\rm (4)} There is a 1--1 correspondence between the isomorphism types and
the graphs of the uniserial left $\la$-modules.

{\rm (5)} The only graphs of uniserial left $\la$-modules are edge
paths. \qed\endproclaim

Clearly every left serial algebra satisfies condition (5) of Theorem 5.6.
Natural instances of algebras satisfying the equivalent conditions of
Theorem 5.6 are, moreover, the algebras $\la =\opio$ where $\O$ is a tiled
classical order over a DVR with uniformizing parameter $\pi$ (see Section 3).

\head 6. Second intermezzo: Saguaros\endhead

Let $\la = \kgam/I$  be a path algebra having a quiver $\Gamma$ without double
arrows, which means that each uniserial module has a unique mast in
$\kgam$.

\definition{Definition 6.1}
Suppose that $T_1, \ldots , T_m$ is a sequence
of non-zero uniserial left $\la$-modules, and let $p_i$ be the mast of $T_i$,
respectively. A left $\la$-module
$T$ is called a  {\it saguaro}\footnote{The name is that of a cactus
found in the Sonoran Desert.} on $(T_1, \ldots , T_m)$ if 

(i) $T\cong \bigl( \bigoplus _{1\leq i\leq m}T_i \bigr) \big/U$,
where $U\subseteq \bigoplus _{1\leq i\leq m}JT_i$ is generated by a sequence
of elements of  the form
$q_it_i-q'_{i+1}t_{i+1}$, $1\leq i\leq m-1$,
where $t_i \in T_i$ are suitable top elements and  $q_i, q'_i$ are right
subpaths of the masts $p_i$ such that $q_it_i\neq 0$, and
$q'_{i+1}t_{i+1}\neq 0$;  moreover, we require that 

(ii) each $T_j$ embeds canonically in $T$ via
$$T_j@>{\roman{can}}>> \biggl( \bigoplus_{1\leq i\leq m}T_i\biggr) \bigg/
U\cong T.$$
The uniserial modules $T_i$ are called the {\it trunks} of
$T$.

In the sequel, we will identify $T$ with
$\bigl( \bigoplus _{1\leq i\leq m}T_i \bigr) \big/ U$. To avoid
ambiguities, we will denote the canonical images of the
trunks $T_i$ by $\That _i$ and the canonical images of the
top elements $t_i$ by $\that _i$.  Any such sequence $(\that_1,
\ldots , \that_m)$ will be called a {\it
canonical sequence of top elements} for $T$.\enddefinition

Note that saguaros are particularly amenable to graphing,
the shape of their graphs explaining their name. By a slight abuse of
language, we will say that the graph of a saguaro $T$ {\it displays a canonical
sequence of top elements} if the simple summands of $T/JT$ shown in the
uppermost layer of the graph are generated by the terms of such a canonical
sequence. It is clear that  layered and labeled graphs displaying a canonical
sequence of top elements always exist. To give an example, the following is a
graph of a saguaro $T$ over the algebra $\la$ of Example 1.1. Since $\gam$ has
no double arrows in this case, we can omit the labels on the edges without
losing information.

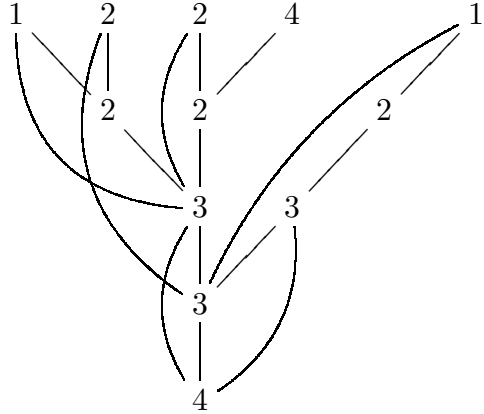
\ignore{
$$\xymatrixcolsep{1.5pc}\xymatrixrowsep{1.5pc}
\xy\xymatrix{
1 \edge[dr] \edge@/_1.75pc/[ddrr] &2 \edge[d] \edge@/_1.75pc/[dddr] &2 \edge[d]
\edge@/_1pc/[dd] &4 \edge[dl] &&1 \edge[dl] \edge@/_1pc/[dddlll]\\
 &2 \edge[dr] &2 \edge[d] &&2 \edge[dl]\\
 &&3 \edge[d] \edge@/_1pc/[dd] &3 \edge[dl] \edge@/^1pc/[ddl]\\
 &&3 \edge[d]\\
 &&4}\endxy$$
}

\noindent Here $T= \bigl( \bigoplus_{i=1}^5 T_i \bigr) \big/U$, where $U$ is
generated by the elements $\alpha t_1-\beta t_2$, $\beta^2t_2-\delta\gamma
t_3$, $\gamma t_3-\gamma t_4$, and $\delta^2\gamma t_4-\delta^2\gamma t_5$.
The trunks of $T$ are the uniserials
$T_1,\dots,T_5$ in $\lamod$ which, relative to suitable top elements
$t_i\in T_i$, have graphs

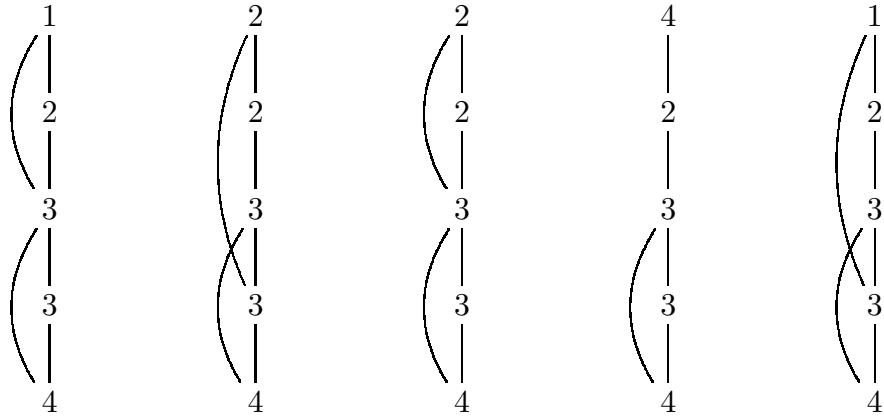
\ignore{
$$\xymatrixrowsep{1.5pc}
\xy\xymatrix{
1 \edge[d] \edge@/_1pc/[dd] &&2 \edge[d] \edge@/_1pc/[ddd] &&2 \edge[d]
\edge@/_1pc/[dd] &&4
\edge[d] &&1 \edge[d] \edge@/_1pc/[ddd]\\
2 \edge[d] &&2 \edge[d] &&2 \edge[d] &&2
\edge[d] &&2 \edge[d]\\
3 \edge[d] \edge@/_1pc/[dd] &&3 \edge[d] \edge@/_1pc/[dd] &&3 \edge[d]
\edge@/_1pc/[dd] &&3 \edge[d] \edge@/_1pc/[dd] &&3 \edge[d] \edge@/_1pc/[dd]\\
3 \edge[d] &&3 \edge[d] &&3 \edge[d] &&3 \edge[d] &&3 \edge[d]\\
4&&4&&4&&4&&4}\endxy$$
}

\noindent respectively.

In Section 8, we will need slight upgrades of some of the observations proved
in \cite{5}; wherever additional care is required, we include the short
proofs for the convenience of the reader.

\proclaim{Observation 6.2}  (On the role of scalars.) If $T\cong
\bigl( \bigoplus_{1\leq i\leq m}T_i \bigr) \big/U$ is
a saguaro, where $U$  is generated by
the relations $q_it_i-q'_{i+1}t_{i+1}$, $1\leq i\leq m-1$, as in
Definition 6.1, and if $k_1, \ldots k_{m-1}$ are non-zero
scalars, then 
$$T\cong \biggl( \bigoplus _{i=1}^m T_i \biggr) \bigg/ \biggl( \sum _{i\leq
m-1}\la (q_it_i-k_iq'_{i+1}t_{i+1}) \biggr). \qed$$ \endproclaim

\proclaim{Observation 6.3}
(On intersections of trunks.) Let
$T=\sum _{1\leq i\leq m}\That _i$ be a saguaro on
$(T_1,\dots,T_m)$ with $q_i\that _i=q'_{i+1}\that _{i+1}$, $1\leq i\leq m-1$,
as in Definition 6.1. Then 
$$\That _i\cap
\That _{i+1}=\la q_i\that _i=\la q'_{i+1}\that _{i+1} \neq 0.$$
In particular, $\soc
(\That _1)=\soc (\That _2)= \dots =\soc (\That_m)\subseteq \soc (T)$.
Moreover, whenever $i<j$, 
$$\That_i\cap \That_j=\That_i\cap \biggl( \sum_{l\geq j}\That_l \biggr) =
\biggl( \sum_{l\leq i}\That_l \biggr) \cap \That_j.$$\endproclaim

\demo{Proof} The first line of equalities is immediate from the
definition. That $\la q_i\that_i \neq 0$ is a consequence of the facts that
$q_it_i$ is nonzero in $T_i$ by condition (i) of Definition 6.1, and that
$T_i$ is isomorphic to its canonical image $\That_i$ in $T$ by condition (ii)
of that definition. This clearly implies that $\soc(\That_i)
=\soc(\That_{i+1})$ for all $i$.

To check that $\That_i\cap \That_j\supseteq \That_i\cap \bigl( \sum_{l\geq
j}\That_l \bigr)$ for $i<j$, let $\lambda \that_i =\sum_{l\geq j}
\lambda_l \that_l$ for suitable elements $\lambda$ and
$\lambda_l$ in $\la$. Note that, whenever $\mu\in\la$ and $k<m$ are such that
$\mu q_kt_k$ is equal to zero in $T_k$, then $\mu q'_{k+1}t_{k+1}$ is zero in
$T_{k+1}$ by condition (ii) of the definition. Using this fact and condition
(i), we obtain $\lambda t_i -\sum_{l\geq j} \lambda_l t_l =\sum_{k\geq i}
\mu_k(q_kt_k -q'_{k+1}t_{k+1})$ in $\oplus_{l\le m} T_l$ for certain
$\mu_k\in\la$. It follows that
$\lambda t_i =\mu_iq_it_i$, and hence that $\lambda\that_i =\mu_iq'_{i+1}
\that_{i+1}\in \That_{i+1}\cap \bigl( \sum_{l\geq j} \That_l \bigr)$, and an
obvious induction on $i$ completes the proof. \qed\enddemo

Note that the inclusion $\soc (\That _i)\subseteq \soc (T)$ is
proper, in general.

\proclaim{Observation 6.4} (Additional information
 on intersections of trunks.) Let $T$ in $\lamod$ be a saguaro on $(T_1,
\ldots , T_m)$ with canonical sequence of top elements $\that_1, \ldots,
\that_m$, where $T_i$ has mast $p_i$. For any two indices $i<j$ in $\{1,\ldots
,m\}$ there then exist right subpaths $a$ of
$p_i$ and $b$ of $p_j$ such that $a\that_i =b\that_j$ and
$\That _i\cap
\That _j=\la a\that_i =\la b\that_j$.\endproclaim

\demo{Proof} We proceed by induction on $j-i$. The case
where $j-i=1$ is covered by Observation 6.3, so we
may assume $j-i\geq 2$. By the induction
hypothesis, we can then find a right subpath $u$ of $p_{i+1}$ and a right
subpath $v$ of $p_j$ with
$u\that_{i+1}=v\that_j$ and such that $\That_{i+1}\cap \That_j=\la
u\that_{i+1}$. In view of Observation 6.3, we see moreover
that 
$$\That_i\cap \That_j=\That_i\cap \That_{i+1}\cap
\That_j =(\That_i\cap \That_{i+1}) \cap(\That_{i+1}\cap \That_j) =\la
q'_{i+1}\that_{i+1}\cap \la u\that_{i+1}$$
where $q'_{i+1}$ is a right subpath of
$p_{i+1}$ as in Definition 6.1.

If $\len(q'_{i+1}) \leq \len(u)$, we obtain $u=wq'_{i+1}$ for a suitable
subpath $w$ of $p_{i+1}$, since both $u$ and $q'_{i+1}$ are right subpaths of
$p_{i+1}$. Observe that, in this case, $wq_i$ is a right subpath of $p_i$ and
$wq_i\that_i =wq'_{i+1}\that_{i+1} =u\that_{i+1} =v\that_j
\neq 0$, which shows in particular that $\That_i\cap
\That_j =\la wq_i\that_i =\la v\that_j$; thus our claim is
satisfied with $a=wq_i$ and $b=v$.

If, on the other hand, $\len(u) <\len(q'_{i+1})$, there exists a subpath $w$
of $p_{i+1}$ with $wu=q'_{i+1}$, which implies $wv\that_j = wu\that_{i+1}
=q'_{i+1}\that_{i+1} =q_i\that_i$ and $\That_i\cap \That_j
=\la q_i\that_i =\la wv\that_j$; in other words, our claim is satisfied with
$a=q_i$ and $b=wv$.
\qed\enddemo

It is an obvious consequence of the preceding observation that, given a
saguaro $T=\sum _{i\in I}\That_i$ over a left serial
algebra and $I_1\subseteq I$, the sum $T'=\sum _{i\in
I_1}\That_i$ is in turn a saguaro with trunks
$\{T_i\mid i\in I_1\}$.

\proclaim{Observation 6.5}  (On the reordering of trunks.) Let $T$
be a saguaro on
$(T_1,
\ldots ,T_m)$, and $\pi \in {\Cal S}_m$ a
permutation. If for each $i\in \{1, \ldots ,m-1\}$ we have
$\That_{\pi (i)}\cap \That _{\pi (i+1)}
\supseteq   \That _{\pi (i)}\cap \That _{\pi (j)}$ for all
$j>i$, then $T$ is also a saguaro on $(T_{\pi (1)}, \ldots
, T_{\pi (m)})$. \endproclaim

\demo{Proof} The details of the proof can be derived from
Observation 6.4 by induction on $m$. \qed\enddemo

In intuitive terms, Observation 6.5 says that, if for all $i<m$
the trunk    $\That _{\pi (i)}$ is linked to $\That _{\pi
(i+1)}$ in at least as high a position in the graph of $T$ as it is linked to
any of the trunks $\That _{\pi (i+1)}, \ldots , \That _{\pi (m)}$,
then $T$  is also a saguaro on
$(T_{\pi (1)}, \ldots ,T_{\pi (m)})$. In particular, this yields: Given any
index $j\in \{1, \ldots ,m\}$, there exists a permutation $\pi \in {\Cal S}_m$
with $\pi (m)=j$ (or with $\pi (1)=j$) such that $T$ is a saguaro on
$(T_{\pi (1)}, \ldots ,T_{\pi (m)})$.

A less obvious consequence of
Observation 6.5 is the following.

\proclaim{Observation 6.6} (Moving two trunks together.)  
Again let $T\in \lamod $ be a saguaro on
$(T_1,\ldots ,T_m)$, and let $(s,t)$ be a pair of distinct indices
in $\{1,\ldots ,m\}$. Then there exists a
permutation $\pi \in {\Cal S}_m$ with the property that $T$ is a
saguaro on $(T_{\pi (1)}, \ldots , T_{\pi (m)})$
and $(s,t)=(\pi (l),\pi (l+1))$ for some $l$. \endproclaim

\demo{Proof} Define $I(s),I(t)\subseteq \{1,\ldots ,m\}$
as follows:
$$\align I(s) &= \{ i \mid i\neq s,t \text{\ and\ } \That_i\cap
\That_s\supsetneq \That_i\cap \That_t\}\\
I(t) &= \{i \mid i\neq s,t \text{\ and\ } \That_i\cap
\That_s\subseteq \That_i\cap \That_t \}.\endalign$$
If $|I(s)|=l-1$ for $l\geq 1$, set $\pi
(l)=s$, and define $\pi (l-1), \ldots ,\pi (1)$
recursively: If $I(s)\neq \varnothing$, i.e., if $l\geq 2$,
select $\pi (l-1)\in I(s)$ so that $\That_{\pi (l-1)}\cap
\That_s$ is maximal among the intersections $\That_i\cap
\That_s$ for $i\in I(s)$. If $I(s)\setminus \{\pi
(l-1)\}\neq \varnothing$, choose an element $\pi (l-2)\in I(s)\setminus \{\pi
(l-1)\}$ so that $\That_{\pi
(l-2)}\cap \That_{\pi (l-1)}$ is maximal among the
intersections $\That _i\cap \That_{\pi (l-1)}$, $i\in I(s)\setminus
\{\pi(l-1)\}$, etc.

Next set $\pi (l+1)=t$, and if $I(t)\neq \varnothing$,
select $\pi (l+2)\in I(t)$ so that $\That_{\pi (l+2)}\cap
\That_t$ is maximal among the $\That_i\cap \That_t, i\in
I(t)$. Continue as above: if $I(t)\setminus \{\pi
(l+2)\}\neq \varnothing$, pick $\pi (l+3)$ in this set
difference so that $\That_{\pi (l+3)}\cap \That_{\pi
(l+2)}$ is maximal among the intersections $\That_i\cap
\That_{\pi (l+2)}$, $i \in I(t)\setminus \{\pi
(l+2)\}$, and so forth.

Going back to the definition of a saguaro, one verifies that $\pi$ satisfies
the hypothesis of Observation 6.5. \qed\enddemo

The category of finite direct sums of saguaros will serve as a key source of
examples to illustrate the usefulness of approximations; see Sections 7 and 8.

\head 7. Second instalment of positive examples: $\pinf(\lamod)$
where $\la$ is left serial\endhead

The importance of saguaros -- or, more generally, modules built on similar
patterns -- came to our attention through their `natural' occurrence as minimal
$\pinf(\lamod)$-ap\-prox\-i\-ma\-tions of the simple modules over a split left
serial algebra $\la$.

\proclaim{Theorem 7.1} {\rm \cite{5, Theorems 5.2, 5.3}} Suppose that
$\la=\kgam/I$ is a left serial algebra. Then $\pinf(\lamod)$ is
contravariantly finite, and the minimal
$\pinf(\lamod)$-ap\-prox\-i\-ma\-tions of the simple left $\la$-modules are
saguaros with simple socles.

More precisely, the minimal
$\pinf(\lamod)$-ap\-prox\-i\-ma\-tion of a simple left $\la$-module $S=\la
e/Je$ can be described as follows: If
$C\subseteq Je$ is maximal with respect to the property that $\pdim \la
e/C<\infty$, there is a unique saguaro $A(S)$ of maximal length in
$\pinf(\lamod)$ such that $\la e/C$ is a trunk of $A(S)$ and $\soc A(S)$ is
simple. Then each canonical epimorphism $A(S)\rightarrow S$, which maps $\la
e/C$ onto
$S$ and sends the other trunks of $A(S)$ to zero, is a minimal
$\pinf(\lamod)$-ap\-prox\-i\-ma\-tion. \qed\endproclaim

When $\la$ is a left serial algebra, the minimal
$\pinf(\lamod)$-ap\-prox\-i\-ma\-tions of the simple modules can actually
be constructed algorithmically from quiver and relations of $\la$. In view of
Theorem 2.4, they form the basic structural components of arbitrary objects
in $\pinf(\lamod)$.

\proclaim{Corollary 7.2} Let $\la$ be a
left serial algebra. Moreover, suppose that the saguaros
$A(S_1),\dots,A(S_n)$ are minimal right $\pinf(\lamod)$-ap\-prox\-i\-ma\-tions
of the simple left $\la$-modules $S_1,\dots ,S_n$. Then a finitely generated
left $\la$-module has finite projective dimension if and only if it is a direct
summand of a submodule $M$ having a filtration $M=M_0\supseteq M_1\supseteq
\cdots\supseteq M_k=0$, where each of the consecutive factors is
isomorphic to some $A(S_i)$. \qed\endproclaim

\example{Example 7.3} Let $\la$ be a left serial algebra whose quiver
has $5$ vertices and whose indecomposable projective modules $\la
e_i$ are given by the following graphs:  

\ignore{
$$\xymatrixcolsep{1.5pc}\xymatrixrowsep{1.5pc}
\xy\xymatrix{
1 \edge[d]^\alpha &2 \edge[d]^\beta &3\edge[d]^\gamma &4 \edge[d]^\delta &5
\edge[d]^\epsilon\\
5 \edge[d]^\epsilon &3 \edge[d]^\gamma &5 \edge[d]^\epsilon &3
\edge[d]^\gamma &5 \edge[d]^\epsilon\\
5 &5 \edge[d]^\epsilon &5 \edge[d]^\epsilon &5 \edge[d]^\epsilon &5\\
 &5 \edge[d]^\epsilon &5 &5 \edge[d]^\epsilon\\
 &5 &&5
}\endxy$$
}

\noindent Then the minimal $\pinf(\lamod)$-ap\-prox\-i\-ma\-tions of the the
simple modules $S_i= \la e_i/Je_i$ are as follows (given the fact that
they are saguaros, they are in fact determined up to isomorphisms by their
graphs): 

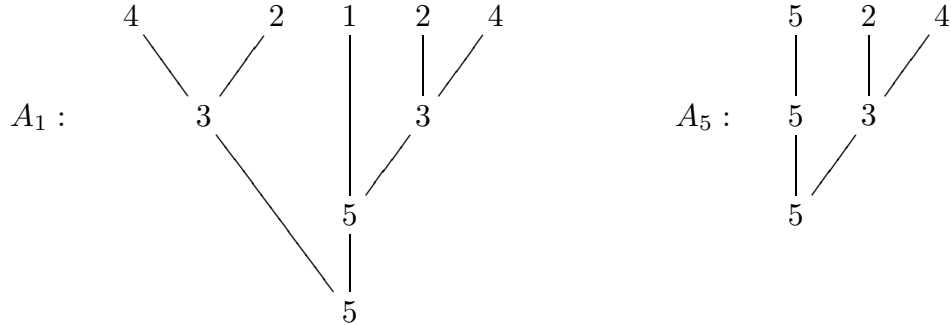
\ignore{
$$\xymatrixcolsep{1pc}\xymatrixrowsep{1.5pc}
\xy\xymatrix{
 &4 \edge[dr] &&2 \edge[dl] &1 \edge[dd] &2 \edge[d] &4 \edge[dl] &&&& 5
\edge[d] &2 \edge[d] &4 \edge[dl]\\
A_1: &&3 \edge[ddrr] &&&3 \edge[dl] &&&&A_5: &5 \edge[d] &3 \edge[dl]\\
 &&&&5 \edge[d] &&&&&&5\\
 &&&&5
}\endxy$$
}

\noindent and $A(S_i)\cong S_i$ for $i=2,3,4$. \qed\endexample

The homological information to be gleaned from saguaros can be pushed
further. For each $d\ge0$, denote by $\Cal P^{(d)}$ the full subcategory of
$\pinf(\lamod)$ whose objects are the modules of projective dimension at most
$d$.

\proclaim{Theorem 7.4} \cite{5} Let $\la$ again be a
left serial algebra. Then, for each $d\geq 0$, the
category ${\Cal P}^{(d)}$ is contravariantly
finite in $\lamod$, and the minimal $\Cal P^{(d)}$-ap\-prox\-i\-ma\-tions of
the simple modules are saguaros. 

A minimal ${\Cal
P}^{(d)}$-ap\-prox\-i\-ma\-tion of the simple  module
$S=\la e/Je$ is as follows: Let $T_d$ be that non-zero homomorphic image of
$\la e$  within
${\Cal P}^{(d)}$ which has smallest (positive) composition length. If $A_d(S)$
is a saguaro in ${\Cal P}^{(d)}$ having simple socle and trunk $T_d$, which has
maximal composition length with respect to these properties, then any
homomorphism
$f : A_d(S)\rightarrow S$ which maps $T_d$ onto $S$ and all the other trunks of
$A_d(S)$ onto zero is a minimal
${\Cal P}^{(d)}$-ap\-prox\-i\-ma\-tion of $S$. In particular, $A_d(S)$ is
unique up to isomorphism.
\qed
\endproclaim

For each simple module $S$, there is thus a sequence of minimal ${\Cal
P}^{(d)}$-ap\-prox\-i\-ma\-tions $A_1(S)$, $A_2(S)$, $A_3(S), \dots$ which
terminates in $A_{\delta}(S)$, where $\delta$ is the left
finitistic dimension of $\la$ (note that
this dimension is known to be finite for left serial algebras
\cite{11, Theorem 3}). Clearly $A_\delta(S)$ coincides with the
minimal $\pinf(\lamod)$-ap\-prox\-i\-ma\-tion of $S$.  As
$d$ increases, the trunks
$T_d$ may shrink from the bottom up, while the saguaros $A_d(S)$ will
increasingly ramify on what is left of these trunks. The following example
illustrates this growth pattern; it encodes a great deal of homological
information about $\la$ in a compact form.

\example{Example 7.5} Let $\la$ be a left serial algebra whose
indecomposable projective modules are represented by the following graphs.

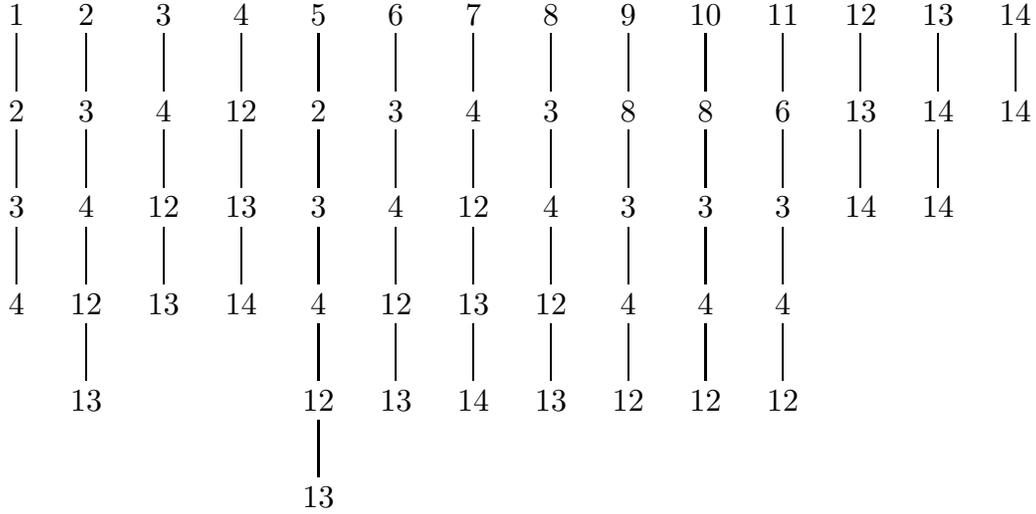
\ignore{
$$\xymatrixcolsep{0.7pc}\xymatrixrowsep{1.5pc}
\xy\xymatrix{
1 \edge[d] &2 \edge[d] &3 \edge[d] &4 \edge[d] &5 \edge[d] &6 \edge[d] &7
\edge[d] &8 \edge[d] &9 \edge[d] &10 \edge[d] &11 \edge[d] &12 \edge[d] &13
\edge[d] &14 \edge[d] \\
2 \edge[d] &3 \edge[d] &4 \edge[d] &12 \edge[d] &2 \edge[d] &3 \edge[d] &4
\edge[d] &3 \edge[d] &8 \edge[d] &8 \edge[d] &6 \edge[d] &13 \edge[d] &14
\edge[d] &14\\
3 \edge[d] &4 \edge[d] &12 \edge[d] &13 \edge[d] &3 \edge[d] &4 \edge[d] &12
\edge[d] &4 \edge[d] &3 \edge[d] &3 \edge[d] &3 \edge[d] &14 &14\\
4 &12 \edge[d] &13 &14 &4 \edge[d] &12 \edge[d] &13 \edge[d] &12 \edge[d] &4
\edge[d] &4 \edge[d] &4 \edge[d] \\
 &13 &&&12 \edge[d] &13 &14 &13 &12 &12 &12\\
&&&&13
}\endxy$$
}

\noindent The evolution of the ${\Cal P}^{(d)}$-ap\-prox\-i\-ma\-tions of the
simple left $\la$-module $S_1$ is graphically represented below.  From left
to right, we exhibit the minimal 
${\Cal P}^{(1)}$-, ${\Cal P}^{(2)}$-, ${\Cal P}^{(3)}$-ap\-prox\-i\-ma\-tions
of
$S_1$; the last coincides with the minimal
$\pinf(\lamod)$-ap\-prox\-i\-ma\-tion, since the left finitistic dimension of
$\la$ is $3$.  The ${\Cal P}^{(0)}$-ap\-prox\-i\-ma\-tion is  simply $\la
e_1$ and is omitted from the list.

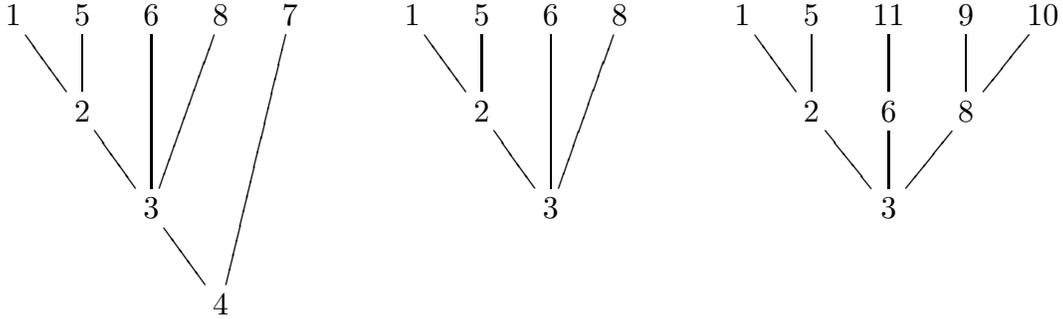
\ignore{
$$\xymatrixcolsep{0.9pc}\xymatrixrowsep{1.5pc}
\xy\xymatrix{
1 \edge[dr] &5 \edge[d] &6 \edge[dd] &8 \edge[ddl] &7 \edge[dddl] &&1
\edge[dr] &5 \edge[d] &6 \edge[dd] &8 \edge[ddl] &&1 \edge[dr] &5 \edge[d]
&11 \edge[d] &9 \edge[d] &10 \edge[dl]\\
 &2 \edge[dr] &&&&&&2 \edge[dr] &&&&&2 \edge[dr] &6 \edge[d] &8
\edge[dl]\\
 &&3 \edge[dr] &&&&&&3 &&&&&3\\
 &&&4
}\endxy$$
}
\endexample

\head 8. Approximations over algebras of finite uniserial type\endhead

Throughout this section we will assume that $\la=\kgam/I$ has finite
uniserial type. In view of Section 5, this implies in particular that $\gam$
has no double arrows. Let $\S\subseteq \lamod$ and $\sinf \subseteq
\lamod$ denote the full subcategories having as objects all finite direct sums
of saguaros in $\lamod$ in the first case, and all finite direct sums of
saguaros of finite projective dimension in the second. If the categories $\S$
and
$\sinf$ have finite representation type, they are of course contravariantly
finite in
$\lamod$ and can thus be used to group the objects of $\lamod$ according to
their minimal
$\S$- or
$\sinf$-ap\-prox\-i\-ma\-tions. Whether this yields an effective classification
depends on the algebra $\la$ and on the class of modules to be explored. The
objects in $\S$ or $\sinf$ are, in a way, `first approximations' to `useful
approximations'. In general, one is led to consider larger contravariantly
finite subcategories of $\lamod$, such as $\Cal C\S\subset \lamod$, the
category of all finite direct sums of saguaros in $\lamod$ and duals of
saguaros in mod-$\la$, or categories of objects glued together along
uniserial submodules from certain cyclic building blocks. An effective method
of creating new, more flexible -- for approximation purposes --
contravariantly finite subcategories from $\S$ is provided by the following
corollary to Theorem 2.5: If $\Cal A\subseteq
\lamod$ is contravariantly finite, then so is the category $\Cal E(\Cal A)$
having as objects all direct summands of extensions of modules in
$\Cal A$ by modules in $\Cal A$.

We believe that the conclusion of the following theorem remains true for
arbitrary algebras of finite uniserial type. The somewhat more narrow
situation which we address here is technically far less involved, however. As
in Section 5, we denote by $V_p$ the affine variety describing the uniserial
left $\la$-modules with mast $p$. For a characterization of the algebras
satisfying the hypothesis of the next theorem, we refer back to
Theorem 5.6.

\proclaim{Theorem 8.1} If $V_p$ is finite for all paths $p$ in $\kgam$, the
subcategories $\S$ and $\sinf$ of $\lamod$ have finite representation
type.\endproclaim

The line of the proof is akin to that of Theorem 5.1 of \cite{BZ-H}, but the
present, far more general, situation calls for non-trivial supplements. We
therefore include a fairly detailed argument.

In case the hypothesis of Theorem 8.1 is satisfied, the graphs of the
saguaros in $\lamod$ are essentially unique, as long as we  insist that a
canonical sequence of top elements be displayed (see Theorem 5.6(4),(5)).  To
state this uniqueness more concisely, we say that two layered and labeled
graphs,
$G_1$ and $G_2$, are {\it equivalent} if  there exists an isomorphism of
undirected graphs between them which preserves the layering, as well as the
numbers attached to the vertices and the labels attached to the edges.  Under
this relation, the equivalence class of graphs  of a saguaro $T$ over an
algebra $\la$ as above is uniquely determined, and  it therefore makes  sense
to refer to {\it the} graph of $T$.

\proclaim{Lemma 8.2} If $V_p$ is finite for all paths $p$ in $\kgam$, then
any two saguaros in $\lamod$ having equivalent graphs (in the sense of
Section 7) are isomorphic. More precisely, if $T^{(1)}$ and $T^{(2)}$ are
saguaros with identical graphs and $\that_{11},\dots,\that_{1s}$
(respectively, $\that_{21},\dots,\that_{2s}$) are canonical sequences of top
elements corresponding to these graphs, then $\that_{1i}\mapsto \that_{2i}$
induces an isomorphism $T^{(1)}\rightarrow T^{(2)}$.\endproclaim

\demo{Proof} Recall that $\la=\kgam/I$ satisfies `$|V_p|<\infty$ for all $p$'
precisely when each graph of a uniserial module is an edge path. In
particular, this means that, given a uniserial module $U\in\lamod$ with top
elements $x$ and $x'$, the left annihilator of $x$ in $\la$ coincides with
that of $x'$; in other words, the assignment $x\mapsto x'$ determines an
automorphism $U\rightarrow U$.

Now suppose that $T^{(1)}$ and $T^{(2)}$ in $\lamod$ are saguaros with
equivalent graphs; without loss of generality, we may assume these graphs to
be identical. Then $T^{(1)}$ and $T^{(2)}$ have identical sequences of
(isomorphism types) of trunks, $(T_1,\dots,T_s)$ say. Write $T^{(j)}=
\sum_{i=1}^n \That_{ji}$ for $j=1,2$, where the $\That_{ji}$ are the
canonical images of the $T_i$ in $T^{(j)}$. Moreover, if $p_i\in\kgam$ is
the mast of $T_i$, the fact that the graphs of $T^{(1)}$ and $T^{(2)}$ are
the same guarantees the existence of right subpaths $q_1$ of $p_1$, $q_i$,
$q'_i$ of
$p_i$ for $2\le i\le s-1$, $q'_s$ of $p_s$ such that $q_i\that_{ji}
=q'_{i+1}\that_{j,i+1}$, $j=1,2$, for canonical sequences of top elements
$\that_{ji}$ corresponding to the coinciding graphs of $T^{(1)}$ and
$T^{(2)}$. By the first paragraph, $\that_{1i}\mapsto
\that_{2i}$ then yields a well-defined isomorphism from $T^{(1)}$ to
$T^{(2)}$. \qed\enddemo

For the sake of the proof of Theorem 8.1 we require some additional
concepts. The reader solely interested in an overview can safely skip them.

\definition{Definition 8.3} Let
$T= \sum_{i=1}^m \That_i$ be a saguaro as in Definition 6.1, and set $I=\{1,
\ldots , m\}$. 

(a) A submodule $V$ of $T$ of the
form  $V=\sum_{i\in I_1}\That _i$, where $I_1\subseteq I$, is
called a {\it subsaguaro} of $T$ if,
for all indices $j\in I\setminus I_1$,
$$\That_j\cap V\subseteq \bigcap_{i\in I_1}\That_i\,.$$

(b) Two subsaguaros $V^{(1)}=\sum _{i\in I_1}\That _i$ and
$V^{(2)}=\sum _{i\in I_2}\That _i$  of $T$ are said to be {\it
isomorphic as subsaguaros} of $T$ if the following is true:  The index sets
$I_1$ and $I_2$ have the same cardinality, $s$ say,  and there exists an
isomorphism $\phi :V^{(1)} \rightarrow V^{(2)}$, together with  orderings
$\{ \That_{11}, \ldots , \That_{1s}\}$ and $\{ \That_{21},\ldots
,\That_{2s}\}$ of the trunks indexed by $I_1$ and by $I_2$, respectively, such
that
$\phi$ induces an isomorphism  $\That _{1j} \rightarrow \That _{2j}$ for $1
\leq j
\leq s$ which restricts to the  identity on the intersection $\bigcap _{i\in
I_1\cup I_2}\That _i$.

(c) The saguaro $T$ is
called  {\it redundant} if it contains nontrivial
isomorphic  subsaguaros $V^{(1)}=\sum _{i\in I_1}\That _i$ and
$V^{(2)}=\sum _{i\in I_2}\That _i$ with $I_1\cap I_2=\varnothing$ such
that, moreover, $\That_k\cap V^{(1)}=\That_k\cap V^{(2)}$ for all $k\in
I\setminus (I_1\cup I_2)$. Otherwise $T$ is called {\it
irredundant}.\enddefinition

Observe that each subsaguaro of $T$ is a saguaro in its own right by the
remark following Observation 6.4. Intuitively speaking, a subsum $V=\sum
_{i\in I_1}\That _i$ of trunks of a saguaro $T=\sum _{i\in I}\That _i$ is a
subsaguaro if there is no index $j\in I\setminus I_1$ such
that the trunk $\That _j$ meets $V$
at a properly higher point in the graph of $T$ than it meets the
intersection $\bigcap_{i\in I_1}\That_i$ of the trunks of $V$.

\proclaim{Lemma 8.4} Again suppose that all the uniserial varieties $V_p$ are
finite, and let $T$ in $\lamod$ be a saguaro with graph $G$. Then $T$ is
redundant if and only if $G$ has a non-trivial full subgraph of the form

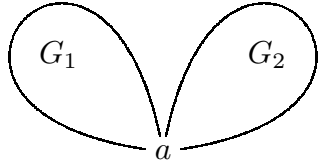
\ignore{
$$\xymatrixcolsep{1.5pc}\xymatrixrowsep{1.5pc}
\xy\xymatrix{
G_1 &&G_2\\
&a \arrow@{-}@'{@+{[0,0]+(-18,2)} @+{[0,0]+(-20,20)} @+{[0,0]+(-3,17)}}
\arrow@{-}@'{@+{[0,0]+(18,2)} @+{[0,0]+(20,20)} @+{[0,0]+(3,17)}}
}\endxy$$
}

\noindent where $G_1$ and $G_2$ are equivalent trees sharing precisely one
vertex $a$ of $G$ such that, moreover, each edge of $G$ which is contiguous
with $G_1\cup G_2$, without belonging to $G_1\cup G_2$, is incident with $a$.
\endproclaim

\demo{Proof} Start by observing that, by Theorem 5.6, our hypothesis forces all
graphs of saguaros in $\lamod$ to be trees (all graphs of uniserial modules
being edge paths). Hence redundancy of $T$ clearly implies the existence of a
full subgraph of $G$ as postulated.

Conversely, suppose that $G$ has a subgraph as described in the claim.
Without loss of generality, we may assume that $G_1$ and $G_2$ are identical
trees. By Observations 6.5 and 6.6, we are, moreover, free to assume that the
graph
$G$ has a form as follows

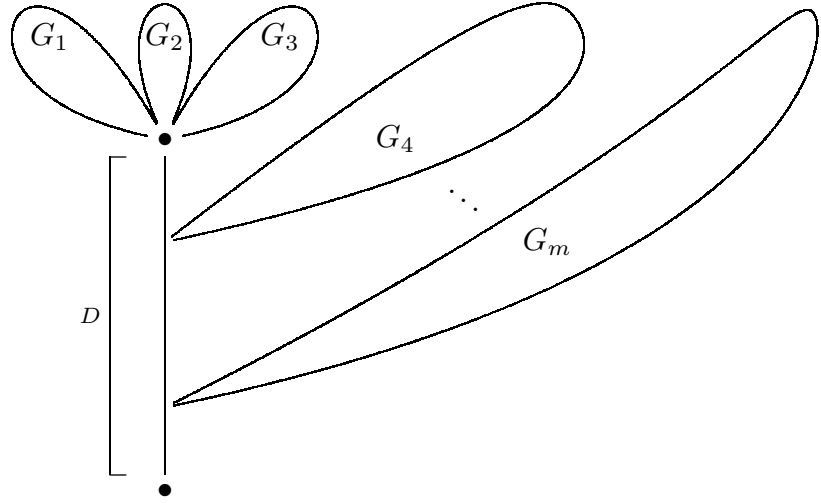
\ignore{
$$\xymatrixcolsep{1.5pc}\xymatrixrowsep{1.5pc}
\xy\xymatrix{
G_1 &G_2 &G_3\\
 &\bullet \arrow@{-}@'{@+{[0,0]+(-17,3)} @+{[0,0]+(-20,20)}
@+{[0,0]+(-7,13)}}
\arrow@{-}@'{@+{[0,0]+(17,3)} @+{[0,0]+(20,20)} @+{[0,0]+(7,13)}}
\arrow@{-}@'{@+{[0,0]+(-5,10)} @+{[0,0]+(0,20)} @+{[0,0]+(5,10)}}
\edge[dddd] \arrow@^{|-|}[dddd]<-4ex>_D &&G_4 \arrow@{}[dr]|\ddots\\
 & \arrow@{-}@'{@+{[0,0]+(47,10)} @+{[0,0]+(50,31)} @+{[0,0]+(35,27)}}
&&&\quad G_m\\
 \\
 & \arrow@{-}@'{@+{[0,0]+(68,14)} @+{[0,0]+(80,65)} @+{[0,0]+(60,30)}}\\
 &\bullet 
}\endxy$$
}

\noindent with a canonical sequence of top elements $\that_{11}, \dots,
\that_{1s}$ corresponding to the uppermost layer of vertices in $G_1$,
 with $\that_{21},
\dots, \that_{2s}$ corresponding to the vertices in the uppermost layer of
$G_2$, and $\that_3,\dots, \that_l$ corresponding to the top vertices of
$G_3,\dots, G_m$. For $j=1,2$, define $\That_{(j,i)} =\la \that_{ji}$ and
$V^{(j)} =\sum_{i=1}^s
\la \That_{(j,i)}$. Then, clearly, $V^{(1)}$ and $V^{(2)}$ are subsaguaros of
$T$, and the sets
$I_j=
\{(j,i)
\mid 1\le i\le s\}$,
$j=1,2$, indexing their trunks are disjoint. We will show that $V^{(1)}$ and
$V^{(2)}$ are isomorphic as subsaguaros of $T$. Indeed, from Lemma 8.2 we
know that the assignment $\that_{1i}\mapsto \that_{2i}$ induces an
isomorphism $\phi : V^{(1)} \rightarrow V^{(2)}$. To see that $\phi$ induces
the identity on $D= \bigcap_{i\le s} \That_{(1,i)} \cap \bigcap_{i\le s}
\That_{(2,i)}$, it suffices to observe that there are right subpaths $a_i$ of
${\operatorname{mast}}(\That_{(1,i)}) ={\operatorname{mast}}(\That_{(2,i)})$,
respectively, such that, for any choice of
$i$ and $h$ in $\{1,\dots,s\}$, we have $\That_{(1,i)}\cap \That_{(2,h)} =D =
\la a_i\that_{1i} =\la a_h\that_{2h}$ and $a_i\that_{1i}
=a_h\that_{2h}$;
for details, see Observation 6.4. \qed\enddemo

\demo{Proof of Theorem 8.1} In a first step we show that, up to isomorphism,
there are only finitely many irredundant saguaros in
$\lamod$. We proceed by induction on the Loewy length $L$
of $\la$. The case $L=1$ being clear, suppose that $L\geq
2$ and that there are $r$ isomorphism types of irredundant
saguaros in $\la /J^{L-1}$-mod; let $G_1,
\dots ,G_r$ be the corresponding graphs, and recall that the $G_i$ are unique
up to equivalence. The graphs of the additional irredundant saguaros in
$\lamod$ are all equivalent to graphs of the form

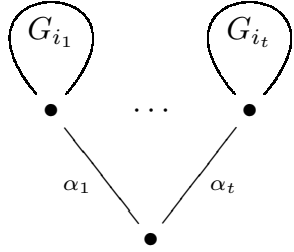
\ignore{
$$\xymatrixcolsep{1pc}\xymatrixrowsep{1pc}
\xy\xymatrix{
G_{i_1} &&G_{i_t}\\
\bullet \arrow@{-}@'{@+{[0,0]+(-8,8)} @+{[0,0]+(0,16)} @+{[0,0]+(8,8)}}
&\cdots &\bullet 
\arrow@{-}@'{@+{[0,0]+(-8,8)} @+{[0,0]+(0,16)} @+{[0,0]+(8,8)}}\\
\\
 &\bullet \edge[uul]^{\alpha_1} \edge[uur]_{\alpha_t}
}\endxy$$
}

\noindent where $t\geq 1$, $1\leq i_j\leq r$, and the $\alpha _i$ are
arrows such that $\alpha _j =\alpha _k$ implies $i_j\ne i_k$. The number of
equivalence classes of graphs of this type is clearly bounded  above by
$2^{ra}$, where $a$ is the number of distinct arrows in $\gam$.

In a second step we prove that each saguaro in $\lamod$
is a direct sum of irredundant ones. Suppose that $T$ is a
saguaro on $(T_1,\dots,T_m)$. This time, we proceed by induction on $m$.
If $T$ is irredundant to begin with, there is nothing to prove; in
particular, this is the case when $m=1$. So suppose that $m\geq
2$ and that  $T$ is redundant with non-zero isomorphic
subsaguaros  $V^{(1)}=\sum _{i\in I_1}\That _i$ and $V^{(2)}=\sum
_{i\in I_2} \That _i$ such that $I_1\cap I_2=\varnothing$ and
$\That_k\cap V^{(1)}=\That_k\cap V^{(2)}\subseteq V^{(1)}\cap V^{(2)}$
for all $k\not\in I_1\cup I_2$ as in  Definition 8.3. 
Moreover, let $\phi :V^{(1)} \rightarrow V^{(2)}$ be an isomorphism and
$\That_{11}, \ldots , \That_{1s}$, resp\. $\That_{21}, \ldots , \That_{2s}$,
orderings of the trunks indexed by $I_1$, resp. by $I_2$, such that $\phi$
induces an isomorphism $\That_{1j} \rightarrow \That_{2j}$ for $ 1 \leq j
\leq s$ and restricts
to the identity on $D:=\bigcap _{i\in I_1} \That _{1i} \cap\bigcap
_{i\in I_2}\That _{2i}$.  In particular, this implies that, given a top
element $\that_{1j} \in \That_{1j}$, the image $\that_{2j} := \phi
(\that_{1j})$ is a top element of $\That_{2j}$ for $1 \leq j \leq s$. 

Observe that $\That_{1i}\cap
\That_{2j}=D$ for all $i,j\in \{1,\ldots ,s\}$, since
$\That_{1i}\cap \That_{2j} \subseteq \bigcap_{1\le k\le s} \That_{1k} \cap
\bigcap_{1\le k\le s} \That_{2k}$ due to the fact that $V^{(1)}$ and $V^{(2)}$
are subsaguaros of $T$ and $I_1\cap I_2 =\varnothing$. It follows that the
submodule
$V:=\sum_{1\leq i\leq s}\la (\that_{1i}-\that_{2i})$ of $T$ is
isomorphic to $V^{(1)}/D$ and thus is a direct sum of saguaros,
each of which has at most $s$ trunks. By the remark following
Observation 6.4, $W:=\sum_{i\in I\setminus
I_1}\That_i$ is in turn a saguaro which clearly has fewer than
$m$ trunks because $I_1\neq \varnothing$. In view  of the
fact that we can reorder $T_1,\ldots ,T_m$ into another
legitimate sequence of trunks in such a way that the trunks
indexed by $I_1$ precede those indexed by $I\setminus I_1$
(Observation 6.5), it is now routine to check that
$V\cap W=0$. We infer that $T=V\oplus W$ and apply the
induction hypothesis to $W$ and to the saguaros occurring as
direct summands of $V$ to complete the proof.  \qed\enddemo

In a first easy example -- Example 4.3 revisited -- we present a monomial
relation algebra
$\la$ for which $\pinf(\lamod)$ fails to be contravariantly finite in
$\lamod$, while all the varieties $V_p$ describing the uniserial left
$\la$-modules are finite. By Theorem 8.1, this implies that $\sinf$ has finite
representation type.

\example{Example 8.5} If $\la$ is the algebra of Example 4.3, then
$\pinf(\lamod)$ fails to be contravariantly finite as we saw earlier. 
However, $\la$ satisfies the hypothesis of Theorem 8.1, and so $\sinf$ has
finite type. A fortiori, $\sinf$ is contravariantly finite in
$\lamod$. The minimal $\sinf$-ap\-prox\-i\-ma\-tions $A_i$ of the simple left
$\la$-modules $S_i$ are as follows:

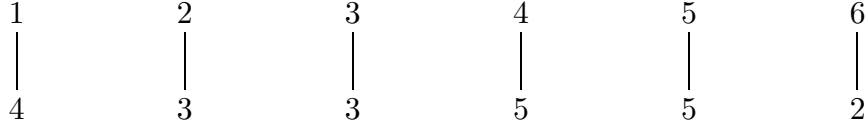
\ignore{
$$\xymatrixcolsep{1.5pc}\xymatrixrowsep{1.5pc}
\xy\xymatrix{
1 \edge[d] &&2 \edge[d] &&3 \edge[d] &&4 \edge[d] &&5 \edge[d] &&6 \edge[d]\\
4 &&3 &&3 &&5 &&5 &&2
}\endxy$$
}

\noindent Observe that $\lfd\la =1=\sup_{1\le i\le 6} \pdim A_i$.
\qed\endexample

\example{Problem 8.6} Characterize those algebras $\la$ of finite uniserial
type for which the supremum of the minimal $\sinf$-ap\-prox\-i\-ma\-tions of
the simple left $\la$-modules equals $\lfd\la$.\endexample

While there are algebras for which this equality fails, for instance among
the algebras of type $\la=\opio$ discussed in Section 4, the
realm of validity of this equality even among non-monomial algebras appears
to be fairly wide. We conclude with a not so straightforward binomial
example. Here the computations leading to the precise shapes of the minimal
approximations of the simple modules are a bit more involved; we will suppress
them nonetheless.

\example{Example 8.7} Let $\la=\kgam/I$, where $\gam$ is the quiver

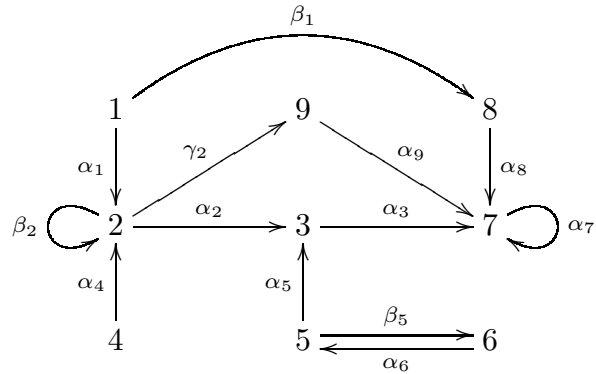
\ignore{
$$\xymatrixcolsep{4pc}
\xy\xymatrix{
1 \arrow[d]_{\alpha_1} \arrow@/^2pc/[rr]^{\beta_1} &9 \arrow[dr]^{\alpha_9}
&8 \arrow[d]^{\alpha_8}\\
2 \lloopd{\beta_2} \arrow[r]^{\alpha_2} \arrow[ur]^{\gamma_2} &3
\arrow[r]^{\alpha_3} &7 \rloopd{\alpha_7}\\
4 \arrow[u]^{\alpha_4} &5 \arrow[u]^{\alpha_5} \arrow[r]<0.5ex>^{\beta_5} &6
\arrow[l]<0.5ex>^{\alpha_6}
}\endxy$$
}

\noindent and the indecomposable projective left $\la$-modules are

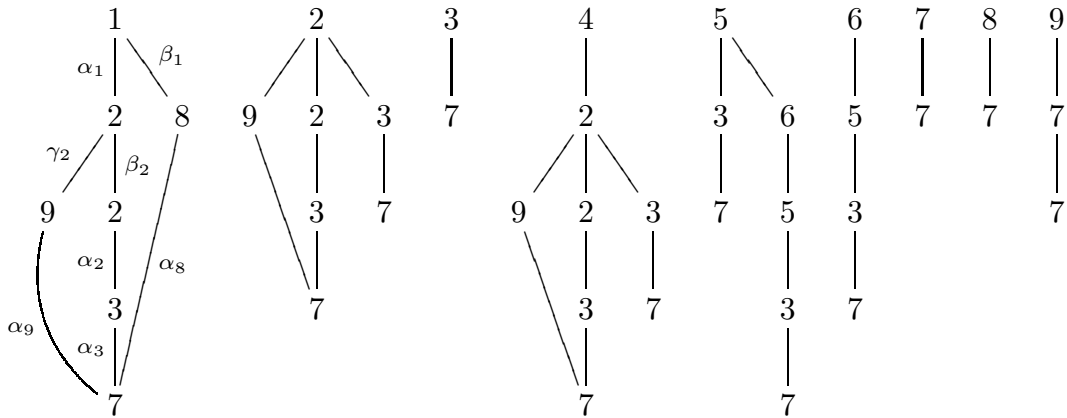
\ignore{
$$\xymatrixcolsep{0.85pc}\xymatrixrowsep{1.5pc}
\xy\xymatrix{
&1 \edge[d]_{\alpha_1} \edge[dr]^{\beta_1} &&&2 \edge[dl] \edge[d]
\edge[dr] &&3 \edge[d] &&4 \edge[d] &&5 \edge[d] \edge[dr] &&6 \edge[d]
&7 \edge[d] &8 \edge[d] &9 \edge[d]\\
&2 \edge[dl]_{\gamma_2} \edge[d]^{\beta_2} &8 \edge[dddl]^{\alpha_8}
&9 \edge[ddr] &2 \edge[d] &3 \edge[d] &7 &&2 \edge[dl] \edge[d] \edge[dr]
&&3 \edge[d] &6 \edge[d] &5 \edge[d] &7 &7 &7 \edge[d]\\
9 \edge@/_1pc/[ddr]_{\alpha_9} &2 \edge[d]_{\alpha_2} &&&3 \edge[d] &7 &&9
\edge[ddr] &2 \edge[d] &3 \edge[d] &7 &5 \edge[d] &3 \edge[d] &&&7\\
&3 \edge[d]_(0.4){\alpha_3} &&&7 &&&&3 \edge[d] &7 &&3 \edge[d] &7\\
&7 &&&&& &&7 &&&7
}\endxy$$
}

\noindent The only relations which cannot be gleaned from the graphs -- the
scalars occurring cannot be detected -- we take to be $\alpha_9\gamma_2
-\alpha_3\alpha_2\beta_2$ and
$\alpha_3\alpha_2\beta_2\alpha_1 -\alpha_8\beta_1$.

Then $\la$ has finite uniserial type, without satisfying the stronger
hypothesis of Theorem 8.1. However, the categories $\S$ and $\sinf$ are
still of finite representation type. One can compute the minimal
$\sinf$-ap\-prox\-i\-ma\-tion
$A_1$ of
$S_1$ to be given by the graph

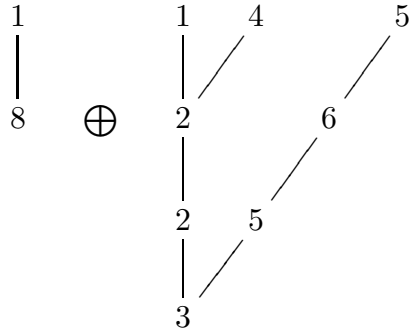
\ignore{
$$\xymatrixcolsep{1pc}\xymatrixrowsep{1.5pc}
\xy\xymatrix{
1 \edge[d] &&1 \edge[d] &4 \edge[dl] &&5 \edge[dl]\\
8 &\bigoplus &2 \edge[d] &&6 \edge[dl]\\
&&2 \edge[d] &5 \edge[dl]\\
&&3
}\endxy$$
}

The minimal $\sinf$-ap\-prox\-i\-ma\-tions $A_i$ of the remaining simples
$S_i$ are:

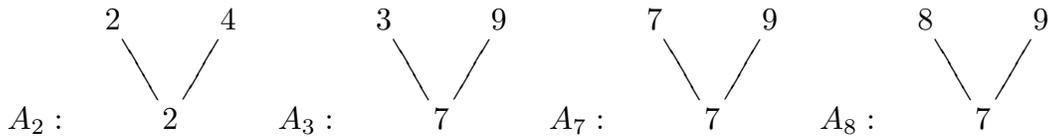
\ignore{
$$\xymatrixcolsep{0.6pc}\xymatrixrowsep{1.5pc}
\xy\xymatrix{
 &2 \edge[dr] &&4 \edge[dl] &&3 \edge[dr] &&9 \edge[dl] &&7 \edge[dr] &&9
\edge[dl] &&8 \edge[dr] &&9 \edge[dl]\\
A_2: &&2 &&A_3: &&7 &&A_7: &&7 &&A_8: &&7
}\endxy$$
}

\noindent and $A_i=S_i$ for $i=4,5,6,9$.

Moreover, one obtains $\sup\{ \pdim A_i \mid 1\le i\le 9\} =2$. We
resolve $A_1$, to indicate how easy it is to find minimal
projective resolutions of saguaros. Indeed, the first syzygy
$\Omega^1(A_1)$ of $A_1$ has graph

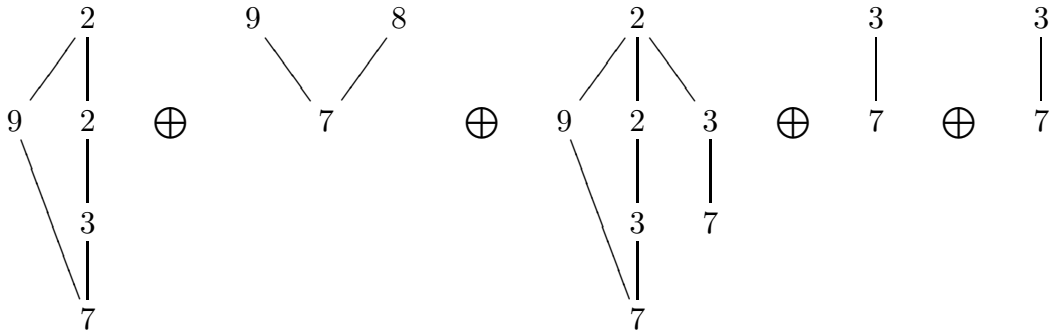
\ignore{
$$\xymatrixcolsep{1pc}\xymatrixrowsep{1.5pc}
\xy\xymatrix{
 &2 \edge[dl] \edge[d] &&9 \edge[dr] &&8 \edge[dl] &&&2 \edge[dl] \edge[d]
\edge[dr] &&&3 \edge[d] &&3 \edge[d]\\
9 \edge[ddr] &2 \edge[d] &\bigoplus &&7 &&\bigoplus &9 \edge[ddr] &2
\edge[d] &3 \edge[d] &\bigoplus &7 &\bigoplus &7\\
 &3 \edge[d] &&&& &&&3 \edge[d] &7\\
 &7 &&&& &&&7
}\endxy$$
}

\noindent and $\Omega^2(A_1)$ has graph

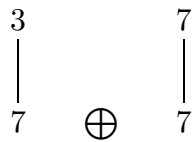
\ignore{
$$\xymatrixcolsep{1pc}\xymatrixrowsep{1.5pc}
\xy\xymatrix{
3 \edge[d] &&7 \edge[d]\\
7 &\bigoplus &7
}\endxy$$
}

\noindent Thus $\Omega^2(A_1)$ is projective. With the `repetition
method' of \cite{8}, one can finally check that $\lfd\la=2$,
which yields the equality $\lfd\la= \sup\{ \pdim A_i \mid
1\le i\le 9\}$. \qed\endexample

\head 9. Phantoms\endhead

When $\A\subseteq \lamod$ fails to be contravariantly finite, we abandon the
requirement that the approximating objects, used to compare arbitrary
finitely generated modules with the modules in $\A$, be themselves finitely
generated. In \cite{10} it became apparent that this often yields
information which is no less effective than that stored in classical (finite
dimensional) $\A$-ap\-prox\-i\-ma\-tions. Such generalized
$\A$-ap\-prox\-i\-ma\-tions, `phantoms' as we will call them, of a simple
module $\la e/Je$ say, provide a synopsis of the relations present in those
objects of $\A$ which contain a top element of type $e$. We simply renounce
the requirement that this picture should fit into a finitely generated module.

\definition{Definitions 9.1} \cite{10} Let $\C\subseteq\A$ be full
subcategories of $\lamod$ and suppose that $\A$ is closed under finite direct
sums. Moreover, let $X$ be a finitely generated left $\la$-module.

(1) A {\it $\C$-ap\-prox\-i\-ma\-tion of $X$ inside $\A$} is a homomorphism $f :
A\rightarrow X$ with $A\in \A$ such that each map $g\in \Hom_\la(C,X)$ with
$C\in\C$ factors through $f$.

(2) An {\it $\A$-phantom of $X$ of the first kind} is an
object $H\in\lamod$ (not necessarily in $\A$) with the following property:
There exists a finite nonempty subclass $\C(H)\subseteq \A$ such that for each
$\C(H)$-ap\-prox\-i\-ma\-tion $f : A\rightarrow X$ inside $\A$, the
module $H$ is a subfactor of $A$. Any direct limit of $\A$-phantoms of $X$ of
the first kind will be called an {\it $\A$-phantom of $X$ of the second kind}.

We will refer to both kinds of phantoms as $\A$-phantoms of $X$.

(3) An $\A$-phantom $H$ of $X$ is called {\it $\C$-effective} if $H$ is a
direct limit of objects in $\A$ and there exists a homomorphism $f :
H\rightarrow X$ with the property that each homomorphism $g\in \Hom_\la(C,X)$
with $C\in\C$ factors through $f$.\enddefinition

The effective phantoms are in a sense the best possible substitutes for
minimal $\A$-ap\-prox\-i\-ma\-tions in the traditional sense. The crucial fact is
that, given any module $X\in\lamod$, nontrivial phantoms  exist. In case
$X$ has an $\A$-ap\-prox\-i\-ma\-tion in the sense of Auslander and Smal\o, the
minimal $\A$-ap\-prox\-i\-ma\-tion $A(X)$ is the `best possible' phantom;
indeed, the
$\A$-phantoms of $X$ are precisely the subfactors of $A(X)$ in that
situation. The case of interest is addressed by the following theorem, the
proof of which indicates a construction pattern for phantoms of infinite
$K$-dimension.

\proclaim{Theorem 9.2} \cite{10} Suppose that $\A\subseteq \lamod$ is
closed under finite direct sums and let $X\in\lamod$. Then the following
conditions are equivalent:

{\rm (1)} $X$ fails to have an $\A$-ap\-prox\-i\-ma\-tion.

{\rm (2)} $X$ has $\A$-phantoms of countably infinite $K$-dimension.

{\rm (3)} There exists a countable subclass $\C\subseteq \A$ such that $X$
has $\C$-effective $\A$-phantoms of infinite $K$-dimension. \qed\endproclaim

\example{Remarks 9.3} (1) From the proof of Criterion 4.1 it can be gleaned
that, under the hypotheses of 4.1, with $p_1,\dots,p_m,q_1,\dots,q_m$ being
paths in
$\kgam$ of positive length, the simple module $S_1$ has an $\A$-phantom with
graph

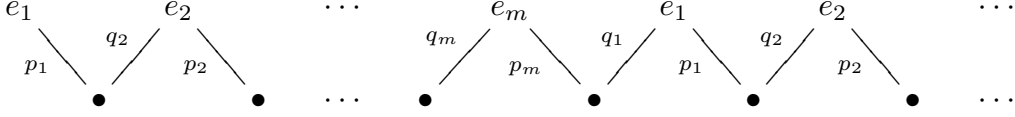
\ignore{
$$\xymatrixcolsep{1pc}\xymatrixrowsep{1.5pc}
\xy\xymatrix{
e_1 \edge[dr]_{p_1} &&e_2 \edge[dl]_{q_2} \edge[dr]_{p_2}
&&\cdots &&e_m \edge[dl]_{q_m} \edge[dr]_{p_m} &&e_1
\edge[dl]_{q_1} \edge[dr]_{p_1} &&e_2 \edge[dl]_{q_2}
\edge[dr]_{p_2} &&\cdots \\
 &\bullet &&\bullet &\cdots &\bullet &&\bullet &&\bullet &&\bullet &\cdots
}\endxy$$
}

The next remark gives a clue how to start building phantoms.

(2) Let $\A=\pinf(\lamod)$ where $\la=\kgam/I$ is a monomial relation
algebra, and suppose that the simple module $S= \la e/Je$ has infinite
projective dimension. If $\alpha_1,\dots, \alpha_m$ are arrows $\alpha_j :
e\rightarrow e_j$ ending in distinct vertices $e_1,\dots, e_m$ such that
$\pdim Je_j =\infty$ for $1\le j\le m$, then

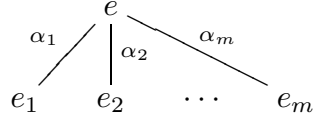
\ignore{
$$\xymatrixcolsep{1pc}\xymatrixrowsep{1.5pc}
\xy\xymatrix{
 &e \edge[dl]_{\alpha_1} \edge[d]^{\alpha_2} \edge[drr]^{\alpha_m}\\
e_1 &e_2 &\cdots &e_m
}\endxy$$
}

\noindent is the graph of a $\pinf(\lamod)$-phantom of $S$.

Moreover, if there exists a module $M\in \pinf(\lamod)$ with graph $G$ having
top elements $m= em$ and $m'= e'm'$ which correspond to the top vertices of
a subgraph of $G$ as follows

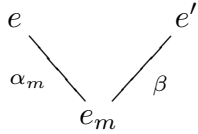
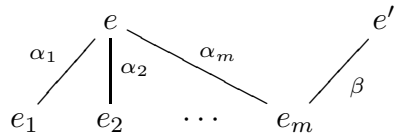
\ignore{
$$\xymatrixcolsep{1pc}\xymatrixrowsep{1.5pc}
\xy\xymatrix{
e \edge[dr]_{\alpha_m} &&e' \edge[dl]^\beta\\
 &e_m}\endxy$$
}

\noindent then

\ignore{
$$\xymatrixcolsep{1pc}\xymatrixrowsep{1.5pc}
\xy\xymatrix{
 &e \edge[dl]_{\alpha_1} \edge[d]^{\alpha_2} \edge[drr]^{\alpha_m} &&&e'
\edge[dl]^\beta\\ 
e_1 &e_2 &\cdots &e_m
}\endxy$$
}

\noindent is the graph of a $\pinf(\lamod)$-phantom of $S$. \qed\endexample

\example{Examples 9.4} We revisit the negative examples of Section 4.
Throughout, $\A$ stands for the category $\pinf(\lamod)$ and $S_1$ for the
simple left $\la$-module corresponding to the vertex `1'.
\medskip

$\bullet$ Example 4.2. The module $M= \varinjlim M_n \in
\pinf(\la\operatorname{-Mod})$ with graph

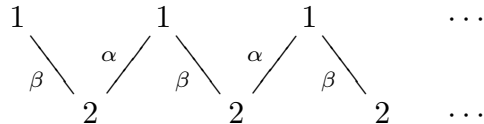
\ignore{
$$\xymatrixcolsep{1pc}\xymatrixrowsep{1.5pc}
\xy\xymatrix{
1 \edge[dr]_\beta &&1 \edge[dl]_\alpha \edge[dr]_\beta &&1 \edge[dl]_\alpha 
\edge[dr]_\beta &&\cdots\\
 &2 &&2 &&2 &\cdots
}\endxy$$
}

\noindent is an $\A$-phantom of $S_1$ which is $\C_0$-effective, where $\C_0=
\{M_n
\mid n\in\N\}$; here $M_n$ is defined as in 4.2.

On the other hand, the left $\la$-module $N$ with graph

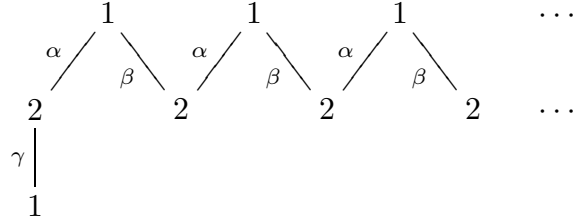
\ignore{
$$\xymatrixcolsep{1pc}\xymatrixrowsep{1.5pc}
\xy\xymatrix{
 &1 \edge[dl]_\alpha \edge[dr]_\beta &&1 \edge[dl]_\alpha \edge[dr]_\beta &&1
\edge[dl]_\alpha \edge[dr]_\beta &&\cdots\\
2 \edge[d]_\gamma &&2 &&2 &&2 &\cdots\\
1}\endxy$$
}

\noindent -- while still an object of $\pinf(\la\operatorname{-Mod})$ with
the property that each homomorphism $M_n\rightarrow S_1$, $n\in\N$, factors
through it -- is not an $\A$-phantom of $S_1$.

Further non-finitely generated phantoms of $S_1$ are as follows. For each
nonzero scalar
$k\in K$ and $n\in\N$, let $x_1 =\dots =x_n =e_1$, and consider the module
$$L_{nk} = \biggl( \bigoplus_{i=1}^n \la x_i \biggr) \biggm/ \biggl( \la
(\beta x_1 -k\alpha x_1) + \sum_{i=2}^n \la (\beta x_i -k\alpha x_i
-\alpha x_{i-1})
\biggr)$$
in $\pinf(\lamod)$, which has the following graph relative to the top
elements $\xbar_1, \dots, \xbar_n$:

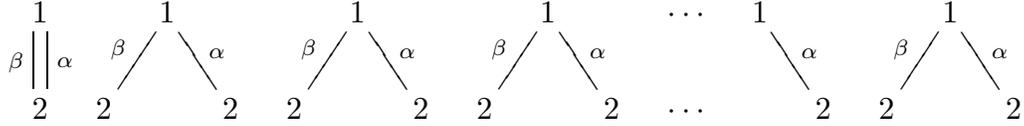
\ignore{
$$\xymatrixcolsep{0.75pc}\xymatrixrowsep{1.5pc}
\xy\xymatrix{
1 \edge[d]<-0.5ex>_\beta \edge[d]<0.5ex>^\alpha &&1 \edge[dl]_\beta
\edge[dr]^\alpha &&&1 \edge[dl]_\beta \edge[dr]^\alpha &&&1 \edge[dl]_\beta
\edge[dr]^\alpha &&\cdots &1 \edge[dr]^\alpha &&&1 \edge[dl]_\beta
\edge[dr]^\alpha\\
2 &2 &&2 &2 &&2 &2 &&2 &\cdots &&2 &2 &&2 
}\endxy$$
}

\noindent Clearly, $L_{nk}$ embeds canonically into $L_{mk}$ for $n<m$,
and it is not difficult to see that the direct limit $L_k :=
\varinjlim L_{nk}$ with graph

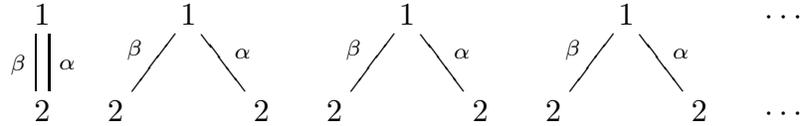
\ignore{
$$\xymatrixcolsep{1pc}\xymatrixrowsep{1.5pc}
\xy\xymatrix{
1 \edge[d]<-0.5ex>_\beta \edge[d]<0.5ex>^\alpha &&1 \edge[dl]_\beta
\edge[dr]^\alpha &&&1 \edge[dl]_\beta \edge[dr]^\alpha &&&1 \edge[dl]_\beta
\edge[dr]^\alpha &&\cdots\\
2 &2 &&2 &2 &&2 &2 &&2 &\cdots 
}\endxy$$
}

\noindent is an $\A$-phantom of $S_1$ which is $\C_k$-effective, where
$\C_k=
\{ L_{nk} \mid n\in\N \}\subseteq \A$.
\medskip

$\bullet$ Example 4.3. Here the direct limit $M= \varinjlim M_n$
with graph

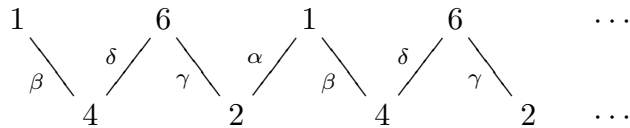
\ignore{
$$\xymatrixcolsep{1pc}\xymatrixrowsep{1.5pc}
\xy\xymatrix{
1 \edge[dr]_\beta &&6 \edge[dl]_\delta \edge[dr]_\gamma &&1 \edge[dl]_\alpha 
\edge[dr]_\beta &&6 \edge[dl]_\delta \edge[dr]_\gamma &&\cdots\\
 &4 &&2 &&4 &&2 &\cdots
}\endxy$$
}

\noindent is an $\A$-phantom of $S_1$ which is $\{ M_n \mid
n\in\N\}$-effective. It is not difficult to
see that the $\A$-phantom $M$ is actually $\A$-effective and hence encodes the
full information stored in classical approximations when they exist.
\medskip 

$\bullet$ Example 4.5. Consider the subclasses $\C= \{M_n \mid n\in\N\}$ and
$\D= \{N_n \mid n\in\N\}$ of $\A$, where $M_n$ and $N_n$ are the left
$\la$-modules with graphs

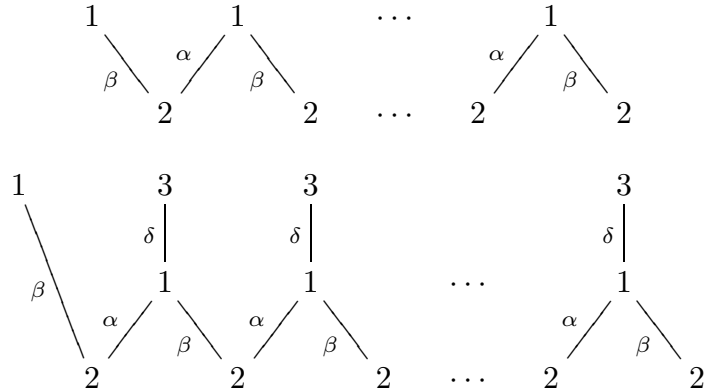
\ignore{
$$\xymatrixcolsep{1pc}\xymatrixrowsep{1.5pc}
\xy\xymatrix{
1 \edge[dr]_\beta &&1 \edge[dl]_\alpha \edge[dr]_\beta &&\cdots &&1
\edge[dl]_\alpha \edge[dr]_\beta \\
 &2 &&2 &\cdots &2 &&2
}\endxy$$

$$\xymatrixcolsep{1pc}\xymatrixrowsep{1.5pc}
\xy\xymatrix{
1 \edge[ddr]_\beta &&3 \edge[d]_\delta &&3 \edge[d]_\delta
&&&&3 \edge[d]_\delta\\
&&1 \edge[dl]_\alpha
\edge[dr]_\beta &&1 \edge[dl]_\alpha \edge[dr]_\beta &&\cdots &&1
\edge[dl]_\alpha \edge[dr]_\beta\\
&2 &&2 &&2 &\cdots &2 &&2
}\endxy$$
}

\noindent respectively, relative to $n$ top elements which are linearly
independent modulo the radical in each case. Then both $M=
\varinjlim M_n$ and $N= \varinjlim N_n$ are $\A$-phantoms
of $S_1$. The phantom $M$ is $\C$-effective, but not $\D$-effective, while
$N$ is $(\C\cup \D)$-effective; in other words, $N$ is the `better' of the
two $\A$-phantoms of $S_1$, storing more information about $\A$ than the
phantom $M$. \qed\endexample

\enddocument